\documentclass{amsart}
\usepackage[english]{babel}
\usepackage{graphicx}
\usepackage[hidelinks]{hyperref} 
\usepackage{amssymb}
\usepackage{amsmath, bm}
\usepackage[foot]{amsaddr}
\usepackage{amsfonts}
\usepackage{mathtools}
\usepackage{algpseudocode}
\usepackage{algorithm}
\usepackage[top=1in, bottom=1.25in, left=1.25in, right=1.25in]{geometry}
\usepackage[table,xcdraw]{xcolor}
\usepackage{multirow, multicol}
\usepackage{tabularx,ragged2e,booktabs}
\usepackage{caption}
\usepackage{subcaption}
\usepackage{hhline}
\usepackage{pgf}
\usepackage{etoolbox}
\usepackage{pgfmath}
\usepackage{enumitem}
\usepackage[sort&compress,numbers]{natbib}
\usepackage{threeparttable}
\usepackage{siunitx}
\usepackage[export]{adjustbox}
\usepackage{scalerel}
\usepackage{tikz}
\usetikzlibrary{external}
\usetikzlibrary{decorations.pathreplacing,calc}

\newcommand{\inmesh}[1]{\scaleto{\mathcal{M}_{#1} \mathstrut}{4pt}}
\newcommand{\inmeshlg}[1]{\scaleto{\mathcal{M}_{#1} \mathstrut}{6pt}}

\newcommand{\inmeshlgtwo}[2]{\scaleto{\mathcal{M}_{#1}^{#2} \mathstrut}{7pt}}
\newcommand\scslash{\stretchrel*{$/$}{\textsc{e}}}
\newcommand{\eitherarrow}{{\leftarrow}\!{\scslash}\!{\rightarrow}}

\graphicspath{{images/}}
\usepackage{float}

\newcommand{\bea}{\begin{eqnarray}}
\newcommand{\eea}{\end{eqnarray}}
\newcommand{\nbea}{\begin{eqnarray*}}
\newcommand{\neea}{\end{eqnarray*}}

\newcommand{\ba}{\begin{align}}

\newcommand{\be}{\begin{equation}}
\newcommand{\ee}{\end{equation}}
\newcommand{\ea}{\end{align}}

\newcommand{\nba}{\begin{align*}}
\newcommand{\nea}{\end{align*}}

\newcommand{\bfmb}{{\mbox{\boldmath{$b$}}}}

\newcommand{\bfmf}{{\mbox{\boldmath{$f$}}}}

\newcommand{\bfmu}{{\mbox{\boldmath{$u$}}}}
\newcommand{\bfmv}{{\mbox{\boldmath{$v$}}}}

\newcommand{\bfmx}{{\mbox{\boldmath{$x$}}}}

\newcommand{\bfmz}{{\mbox{\boldmath{$z$}}}}
\newcommand{\bfmA}{{\mbox{\boldmath{$A$}}}}

\newcommand{\bfmV}{{\mbox{\boldmath{$V$}}}}
\newcommand{\bfmW}{{\mbox{\boldmath{$W$}}}}

\usepackage{tikz,xcolor}
\definecolor{lime}{HTML}{A6CE39}
\DeclareRobustCommand{\orcidicon}{%
	\begin{tikzpicture}
	\draw[lime, fill=lime] (0,0) 
	circle [radius=0.16] 
	node[white] {{\fontfamily{qag}\selectfont \tiny ID}};
	\draw[white, fill=white] (-0.075,0.095) 
	circle [radius=0.007];
	\end{tikzpicture}
	\hspace{-2mm}
}
\foreach \x in {A, ..., Z}{%
	\expandafter\xdef\csname orcid\x\endcsname{\noexpand\href{https://orcid.org/\csname orcidauthor\x\endcsname}{\noexpand\orcidicon}}
}

\begin{document}
\title[GFN for reduced operator learning in multifidelity applications]{GFN: A graph feedforward network for resolution-invariant reduced operator learning in multifidelity applications}

\author{Oisín M. Morrison$^1$ \orcidA, Federico Pichi$^{1,2}$ \orcidB, and Jan S. Hesthaven$^1$ \orcidC}
\address{$^1$ Chair of Computational Mathematics and Simulation Science, École Polytechnique Fédérale de Lausanne, 1015 Lausanne, Switzerland}
\address{$^2$ mathLab, Mathematics Area, SISSA, via Bonomea 265, I-34136 Trieste, Italy.}

\begin{abstract}
This work presents a novel resolution-invariant model order reduction strategy for multifidelity applications. We base our architecture on a novel neural network layer developed in this work, the graph feedforward network, which extends the concept of feedforward networks to graph-structured data by creating a direct link between the weights of a neural network and the nodes of a mesh, enhancing the interpretability of the network. We exploit the method's capability of training and testing on different mesh sizes in an autoencoder-based reduction strategy for parametrised partial differential equations. We show that this extension comes with provable guarantees on the performance via error bounds. The capabilities of the proposed methodology are tested on three challenging benchmarks, including advection-dominated phenomena and problems with a high-dimensional parameter space. The method results in a more lightweight and highly flexible strategy when compared to state-of-the-art models, while showing excellent generalisation performance in both single fidelity and multifidelity scenarios.

\medskip
\noindent \textbf{Keywords:} \textit{graph neural networks}, \textit{model order reduction}, \textit{operator learning}, \textit{resolution invariance}, \textit{multifidelity surrogate modelling}, \textit{parametrised PDEs}.

\medskip
\begin{center}
    \textbf{Code availability:} \url{https://github.com/Oisin-M/GFN}
\end{center}
\end{abstract}

\maketitle

\section{Introduction and Motivation}
Traditional numerical solvers for the high-fidelity approximation of partial differential equations (PDEs) are computationally prohibitive in real-time and many-query contexts, underpinning the need for quicker evaluations of the numerical solution of PDEs. Reduced order models (ROMs) have arisen as a means to address this issue for parametrised PDEs, accelerating the process via the creation of efficient computational models \cite{Benner2017}, with fewer degrees of freedom.

Model order reduction (MOR) is therefore concerned with the creation of offline-online surrogate models that are cheaper whilst maintaining high levels of accuracy. One such popular approach is the reduced basis method, characterised by the creation of a reduced space, typically obtained via proper orthogonal decomposition (POD) or greedy methods \cite{Quarteroni2015, Hesthaven2015}. Such methods are generally linear, rendering them inefficient for problems exhibiting a slow decay in the Kolmogorov $n$-width \cite{Greif2019, Ohlberger2015}. To combat this, a range of nonlinear ROM techniques have been developed, with a notable subclass of such methods being machine learning-based approaches. In particular, much promise has been shown in incorporating autoencoder-based architectures into the ROM context \cite{Lee2020, Fresca2021, Pichi2023, Romor2023} due to their feasibility as a nonlinear extension of principal component analysis \cite{Kramer1991}, displaying superior compression capabilities compared to linear methods especially when dealing with problems with slow Kolmogorov $n$-width decay \cite{Greif2019, Ohlberger2015}.

Such autoencoder-based MOR strategies typically benefit from being non-intrusive approaches, removing the requirement for any knowledge of the generation procedure for the training data, therefore leading to possible exploitation of experimental data \cite{Forrester2010, Fidkowski2014, Kuya2011}.
Moreover, contrarily to physics-informed neural networks (PINNs) \cite{Raissi2019} or sparse identification of nonlinear dynamics (SINDy) \cite{Brunton2016}, these are usually not physics-based methods, and hence render it possible to learn the solution's behaviour even without prior knowledge.

In this direction, some approaches seek to recover the operator mappings, also known as operator inference, on some fixed discretised domain \cite{Peherstorfer2016, Qian2022, Kramer2024}. Instead of approximating the solution at some fixed locations, advanced neural operators (NOs) go further by directly approximating the solution operator itself \cite{Kovachki2021,BoulleMathematicalGuideOperator2023}. As a result, NOs benefit from the desirable property of being independent from a fixed discretisation, possibly achieving super- and sub-resolution \cite{Li2020c, Raonic2023, Sun2022, Kovachki2021}.

The development of so-called \textit{resolution-invariant} methods, i.e.\ methods without a dependence on a fixed discretisation, can be of the utmost importance when dealing with experimental data. This is especially the case in climate modelling, where information depends on both weather stations and satellites, the number and locations of which can vary greatly \cite{Dinku2018, Dinku2019, Mukkavilli2023}. These considerations have lead to similar and complementary approaches including the aforementioned NOs, graph neural networks (GNNs) and multifidelity reduction strategies.

Multifidelity ROMs are typically focused on leveraging cheaper computational data, known as low-fidelity data, in addition to more expensive high-fidelity data, to build an efficient offline MOR strategy. Low-fidelity may be interpreted as simplification of the governing equations, low-order approximations, or scarce information, but here we focus on the case of coarser discretisations, aiming to obtain a ROM capable of learning from data on multiple fidelities. State-of-the-art approaches for ROM \cite{Guo2022, Kast2020, Meng2020, Conti2023}, however, are not generally capable of learning from any arbitrary resolution, but rather from a prespecified fixed set of possible discretisations.

On the other hand, GNNs are capable of leveraging data from graphs of arbitrary cardinality, while preserving the structure of the problem by enforcing geometric priors into the model. A large number of such works (reviewed in \cite{Zhang2018, Zhang2019, Liu2022}) have been focused on the generalisation of pooling, unpooling and convolutional layers to graphs, achieving excellent results across a variety of domains including drug development \cite{Sun2020}, traffic prediction \cite{Jiang2022}, finance \cite{JWang2021}, climate \cite{Lam2023} and fluid mechanics \cite{Hernandez2022, Belbute2020}.
Moreover, being designed for unstructured data, GNNs have an advantage over many current state-of-the-art NOs such as Fourier neural operators \cite{Li2020c} and convolutional neural operators \cite{Raonic2023} which require data on structured grids.

\begin{figure}[!ht]
    \centering
\includegraphics[width=0.92\textwidth]{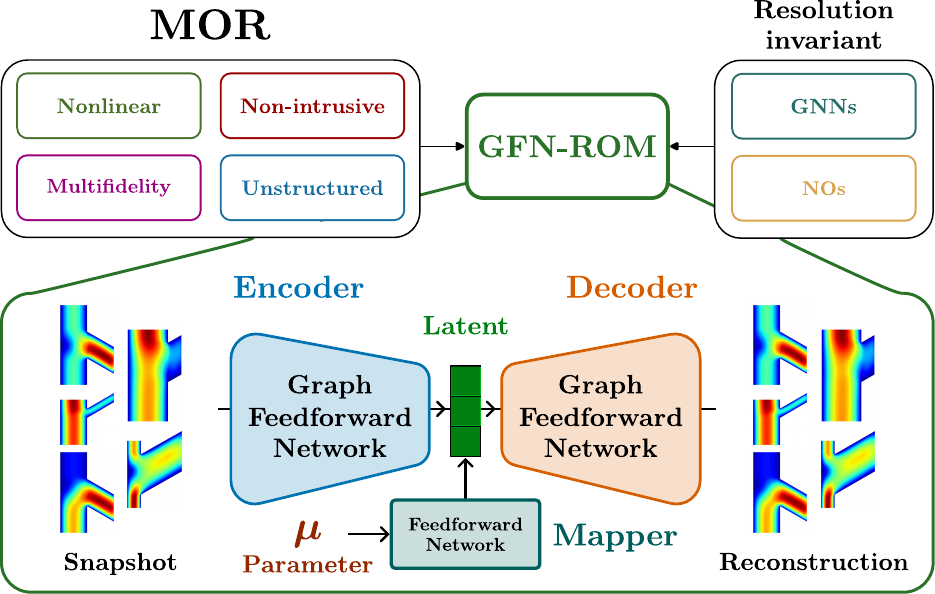}
    \caption{GFN-ROM is a nonlinear non-intrusive multifidelity ROM capable of dealing with unstructured data, interplaying between MOR and resolution-invariant techniques.} 
    \label{fig:gfn_rom_scheme}
\end{figure}

\textbf{Contributions.} In this work, we make the following contributions:
\begin{enumerate}[label=(\roman*)]
    \item We introduce a novel neural network layer called the \textit{graph feedforward network} (GFN) which generalises the feedforward layer for data on arbitrary discretisations.
    \item We introduce \textit{GFN-ROM}, a new autoencoder-based ROM approach for multifidelity data using GFNs.
    \item We provide bounds on the super- and sub-resolution errors using GFNs and GFN-ROM.
    \item We demonstrate the remarkable generalisation properties of GFN-ROM through a collection of challenging problems.
\end{enumerate}

To the best of our knowledge, GFN-ROM is the only graph-based resolution-invariant ROM. By associating weights and biases to nodes in a mesh, GFN-ROM also allows for better interpretability compared to approaches such as DL-ROM or GCA-ROM. Furthermore, we show that GFN-ROM is closely related to the standard definitions of NOs and GNNs, by retaining important properties of both concepts (see \autoref{fig:gfn_rom_scheme}).
    
The remainder of the paper is organised as follows. \hyperref[sec:ROM]{Section~\ref*{sec:ROM}} presents the standard intrusive and non-intrusive framework for MOR. \hyperref[sec:MROM]{Section~\ref*{sec:MROM}} defines the challenges in multifidelity applications, while \hyperref[sec:resinv_method]{Section~\ref*{sec:resinv_method}} introduces resolution invariant methods. The core of the manuscript is presented in \hyperref[sec:GFN]{Section~\ref*{sec:GFN}} where the GFN is introduced to leverage data on any discretisation. \hyperref[sec:GFN-ROM]{Section~\ref*{sec:GFN-ROM}} and \hyperref[sec:adaptive_training]{Section~\ref*{sec:adaptive_training}} present the exploitation of GFN in the ROM context, and how to perform adaptive multifidelity training. Finally, \hyperref[sec:results]{Section~\ref*{sec:results}} shows the performance of the proposed methodology on complex MOR benchmarks, and \hyperref[sec:conc]{Section~\ref*{sec:conc}} gives some conclusions and highlights future perspectives.

\section{Reduced Order Models}\label{sec:ROM}

Full order models are designed to solve PDEs via high-fidelity systems of equations. Such systems involve a large number of degrees of freedom $N_h$, which makes their numerical solution costly to compute. In a parametrised PDE context, this is particularly prohibitive since one aims at recovering solutions in real-time for a number of different physical or geometrical configurations, giving rise to a parameter-dependent high-fidelity solution $\bfmu_h(\boldsymbol{\mu})$. ROMs seek to mitigate this issue by building and subsequently solving a cheaper reduced system with a much lower number of degrees of freedom $N \ll N_h$, giving rise to the reduced order solutions $\bfmu_N(\boldsymbol{\mu})$. High-fidelity solutions can then be approximated from the reduced order ones via some transform $\bfmu_h(\boldsymbol{\mu}) \approx \phi \left( \bfmu_N(\boldsymbol{\mu}) \right)$, where $\phi$ is a linear or nonlinear mapping dependent on the choice of ROM employed. Such methods can be constructed in an intrusive or non-intrusive fashion, depending on whether knowledge of the high-fidelity system is required for the computation of $\bfmu_N(\boldsymbol{\mu})$.

\subsection{Projection-based intrusive methods}

Intrusive methods require knowledge of the high-fidelity system. For ease of illustration, we consider the case of a linear high-fidelity system:
\bea
\label{eq:high_fid}
\bfmA_h(\boldsymbol{\mu}) \bfmu_h(\boldsymbol{\mu}) = \bfmf_h(\boldsymbol{\mu}),
\eea
where $\boldsymbol{\mu} \in \mathbb{R}^{N_\mu}$ is the parameter, $\bfmA_h \in \mathbb{R}^{N_h \times N_h}$ is the stiffness matrix, $\bfmu_h \in \mathbb{R}^{N_h \times  N_{V_F}}$ is the high-fidelity solution, and $\bfmf_h \in \mathbb{R}^{N_h \times N_{V_F}}$ is the forcing term\footnote{For ease of notation, in the following we drop the dependence on the parameters $\boldsymbol{\mu}$.}.

A ROM aims at replacing \autoref{eq:high_fid} with a cheaper reduced system. Linear MOR techniques express an approximation of the high-fidelity solution $\tilde{\bfmu}_h$ as a linear expansion over some chosen basis functions $\{\boldsymbol{\psi}_i\}_{i=1}^N$ with $\boldsymbol{\psi}_i \in \mathbb{R}^{{N_h}}$. Introducing the matrix $\bfmV = \left[ \boldsymbol{\psi}_1 | \cdots | \boldsymbol{\psi}_N \right]$ one can concisely write $\bfmu_h \approx \bfmV \bfmu_N$. POD is a popular choice for constructing a basis, extracting the principal components from a series of high-fidelity solutions for different parameter realisations, known as snapshots. POD provides the best rank $N$ subspace for approximating the snapshots in a least-squares sense \cite{Manzoni2016}. Once the set of basis functions has been selected, it remains to calculate the reduced coefficients $\bfmu_N$. Intrusive projection-based MOR approaches exploit the knowledge of the high-fidelity system by imposing the projected residual of the reduced order solution onto $\bfmV$ is zero i.e.\ one solves for $\bfmu_N$ the system
\bea
\label{eq:rom_force}
\bfmV^T (\bfmf_h - \bfmA_h \bfmV \bfmu_N) = \boldsymbol{0}.
\eea
An efficient application in the many-query context requires \autoref{eq:rom_force} to exhibit an affine parametric dependence in order to circumvent $\mathbb{R}^{N_h}$-dependent online operations. This is usually not fulfilled when dealing with complex and/or nonlinear problems.

\subsection{Non-intrusive methods}
When the high-fidelity system is not available, e.g.\ black-box or commercial solvers, intrusive MOR cannot be pursued and the reduced coefficients have to be recovered in a non-intrusive manner. For example, POD with interpolation (PODI) recovers the reduced coefficients from the projection of the snapshots onto the POD basis \cite{Bui2003, Demo2019}. A similar regression-based  machine learning approach, POD-NN, is also possible by training a network to predict the reduced coefficients, with the projected coefficients used as the training dataset \cite{Hesthaven2018, Barnett2022}.

Although non-intrusive, PODI, POD-NN and other such methods still represent linear approaches since they approximate the high-fidelity coefficients via a linear scheme $\bfmu_h \approx \bfmV \bfmu_N$. This can be problematic, since many important physics phenomena exhibit a slow decay in the Kolmogorov $n$-width, making linear methods inefficient for such cases \cite{Greif2019, Ohlberger2015}. As a result, nonlinear approaches have arisen which seek to avoid this limitation by instead prescribing the more general form of $\bfmu_h(\boldsymbol{\mu}) \approx \psi \left( \bfmu_N(\boldsymbol{\mu}) \right)$, where $\psi$ is a nonlinear mapping. A huge variety of nonlinear compression schemes have been explored for standard MOR: using local reduced bases to construct a piecewise linear scheme \cite{Amsallem2012}, considering nonlinear compression via kernel POD \cite{Diez2021,KhamlichOptimalTransportinspiredDeep2023}, shifted POD \cite{Reiss2018} or registration methods \cite{Taddei2020}.

Machine learning approaches are also possible, with the advantage of typically being cheap at inference. In this direction, Gaussian process regression (GPR) approaches have been heavily investigated \cite{Guo2018, Guo2019, Zhang2019a, Kast2020, Guo2022}. In contrast to standard MOR, these approaches learn a probability distribution and are also advantageous in terms of uncertainty quantification \cite{Cicci2023}. However, GPR models are expensive with cubic scaling, and direct application of GPR to large datasets is infeasible, instead requiring techniques such as local approximate GPR \cite{Gramacy2016, Gramacy2020}.

On the other hand, autoencoders have also arisen as a hugely popular means of nonlinear compression and do not suffer from such poor scaling. The method was originally introduced as a nonlinear generalisation of principal component analysis \cite{Kramer1991}, making the approach a generalisation of POD-based methods. Recently, it has been used in many MOR applications with great success \cite{Moya2022, Hernandez2021, Lee2020, Romor2023, Fresca2021, FrancoDeepLearningbasedSurrogate2023a}, which makes it a hugely promising direction that we pursue in this work.

\section{Multifidelity Reduced Order Models}\label{sec:MROM}

By removing any dependence on the generation procedure for the training data, non-intrusive ROMs can be used as black-box approaches, allowing e.g.\ to augment a small experimental dataset with simulated data \cite{Forrester2010, Fidkowski2014, Kuya2011}. In fact, even without embedding any physics into the ROM, it is still possible to learn the dynamics of the system. Nonetheless, standard approaches typically require all training data to be generated or measured on the same grid, and with the same accuracy. Whilst often assumed for algorithmic convenience, these requirements are highly undesirable for two key reasons: firstly, generating highly-fidelity data for training is expensive and secondly, experimental data is often multifidelity data.

\subsection{Generating high-fidelity data for training is expensive} Several problems require fine discretisations to obtain accurate results and resolve multiscale systems.
However, ROMs need to sample the parameter space well in order to make reasonable predictions. Thus, for standard MOR this involves the expensive computation of a large number of high-fidelity solutions. By allowing a ROM to train on data of different fidelities, one has the potential to greatly reduce the computational burden by computing only a small number of expensive (high-fidelity) simulations, and augmenting the data with a number of cheap (low-fidelity) simulations. This reflection has lead to the development of a large number of multifidelity ROMs \cite{Guo2022, Kast2020, Meng2020, Raissi2017, Conti2023, Fernandez2023, Heiss2023}, aiming to alleviate the computational burden of generating training data.

Low-fidelity simulated data can be generated in a number of ways, either by reducing the accuracy of the computational method e.g.\ using $\mathbb{P}1$ instead of $\mathbb{P}2$ finite elements, considering a coarsened discretisation, or by simplifying the physics. The central task for multifidelity MOR is to determine how to combine data arising from a number of different fidelities. A common assumption is that all fidelities are prespecified and separate models are built for each of the $F$ different fidelities. A common approach has been to adopt a multifidelity GPR based approach for MOR \cite{Raissi2017, Guo2022, Kast2020, Meng2020, Conti2023}, commonly known as cokriging \cite{Kennedy2000, Alvarez2012}, expressing each of the $F$ solution fidelities as a combination of $F$ independent Gaussian processes. Similarly, multigrid approaches have also been heavily explored for parametric PDEs \cite{Teckentrup2015, Harbrecht2016, Ballani2017, Lye2021, Heiss2021, Heiss2023}, relying on a hierarchy of models, and iteratively computing corrections through a sequence of $F$ coarser discretisations.

\subsection{Experimental data is often multifidelity data} Other than synthetic data, ROMs are of utmost importance also when dealing with experimental data, e.g.\ coming from sensors which may move or break. This can lead to a high number of resolutions $F$, particularly for contexts such as climate modelling where data is dependent on weather stations, the number and locations of which can vary greatly \cite{Dinku2018, Dinku2019}. As a concrete example, the average number of stations used in the Global Precipitation Climatology Center \cite{Becker2013} full-data product over Africa has varied from under 250 to over 3250 stations between 1901 and 2013 \cite{Dinku2019}. Even with a coarse yearly temporal resolution, this would still give a total of $F=112$ different resolutions.

This is particularly problematic when one considers multifidelity ROMs train a new model for each fidelity. Furthermore, in many cases it is not possible to upscale or downscale the model predictions to unseen fidelities (e.g.\ finer discretisation), or to leverage training data from other previously unseen fidelities either.

\begin{figure}[!h]
    \centering
    \includegraphics[width=.98\textwidth]{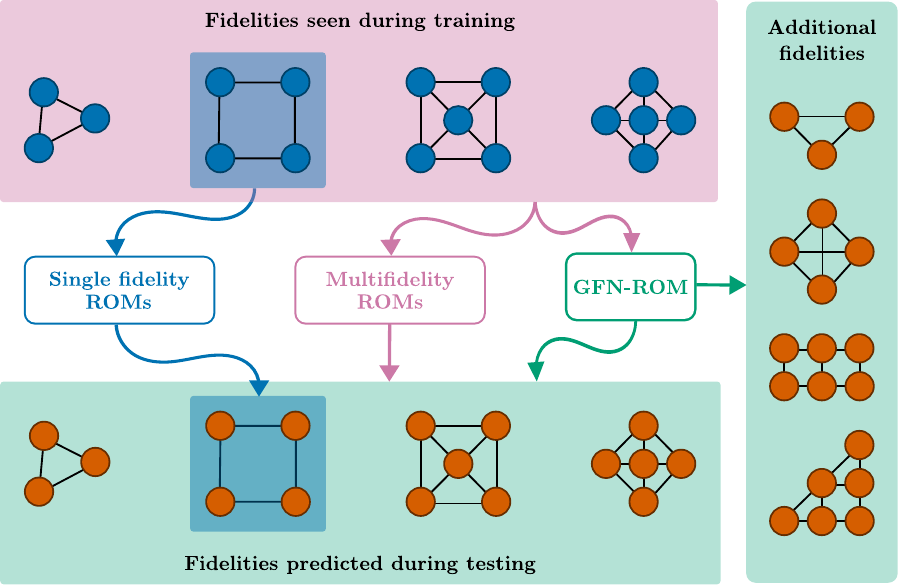}
    \caption{GFN-ROM as a resolution-invariant ROM, capable of handling data from any discretisation, both in training and in testing modes.}
    \label{fig:GFN_res_inv}
\end{figure}

Thus, in the experimental context, a single model capable of training and being evaluated on complexity arbitrary discretisations is required as depicted in \autoref{fig:GFN_res_inv}, i.e.\ there should be no dependence in the model architecture on the discretisations of the training samples. The ability to evaluate and train with data on any discretisation is known as resolution invariance \cite{Kovachki2021}.

None of the MOR approaches discussed thus far are resolution-invariant, underscoring the need for the development of a new, resolution-invariant ROM.

\section{Related Resolution-Invariant Methods}\label{sec:resinv_method}

There exist a number of resolution-invariant approaches relevant for PDEs problems. Of particular interest are PINNs, NOs, and GNNs.

\subsection{Physics-informed neural networks}
Standard ROM approaches usually do not require any explicit positional information, simply relying on the snapshot matrix, thereby making resolution-invariant extensions not straightforward. However, already when considering the approximation of a single snapshot, it is possible to map directly from positional information to a PDE solution i.e.\ $\bfmx \mapsto \phi(\bfmx) \approx u(\bfmx)$, where $\phi$ is a function to fit.

This is precisely the idea of PINNs \cite{Raissi2018, Raissi2019, Cuomo2022}, whereby a neural network $\phi$ is trained via a physics-informed loss i.e.\ the loss is taken to be the residual of the network approximation for the governing PDE, and no training data is required.  Whilst PINNs can perform well on several problems \cite{Raissi2019a, Raissi2019b, Raissi2020, Sun2020a, Jin2021, Arzani2021, DeRyck2022, Kashefi2022}, they suffer from a number of issues such as ``spectral bias”, where they are incapable of learning functions with high frequencies \cite{Wang2022a}, convergence discrepancies among their loss components \cite{Wang2022a} and stiff gradient flow dynamics \cite{Wang2021a}. Furthermore, their extension towards parametrised problems is in practice unsuitable due to the complexity of the training phase, rendering them not useful in the ROM context \cite{Li2021, Rosofsky2023a}. For example, the POD-PINN approach trains a neural network to predict the reduced coefficients from the parameters by minimising the residual of the reduced problem, making it an intrusive approach requiring an efficient means of evaluating the residual of the reduced problem \cite{Chen2021}.

\subsection{Neural operators}\label{subsec:nos}
The inability of PINNs to learn several instances for parametric PDEs is due to the fact that such methods are not designed to directly model the solution operator. Unlike PINNs and standard machine learning methods which learn a map between evaluations of two functions, neural operators learn maps between the function spaces themselves \cite{Kovachki2021}. By attempting to model explicitly the solution operator itself, neural operators remove the burdensome dependence of typical machine learning methods on the discretisation of the data, meaning they are inherently capable of achieving super- or sub-resolution \cite{Li2020c, Raonic2023, Sun2022, Kovachki2021} and also typically benefit from desirable properties such as discretisation invariance \cite{Kovachki2021}. Neural operators have shown impressive performance \cite{Raonic2023, Li2020c, Hao2023, Cao2021, Bonev2023}, even learning different physics at the same time \cite{Kurth2023, Pathak2022}. However, many state-of-the-art neural operators such as Fourier neural operators \cite{Li2020c} and convolutional neural operators \cite{Raonic2023} require structured data as input. Additionally, such methods are designed for learning maps between two spatially dependent functions, and not from a global parameter to a PDE solution, as in the simplified yet more common MOR setting. With the exception of the general neural operator transformer \cite{Hao2023}, existing neural operators therefore do not take global parameters as input.

\subsection{Graph neural networks}\label{subsec:gnns}
PDEs are often posed on complex domains, giving rise to unstructured meshes when computing numerical solutions. Similarly, experimental data is also rarely present on structured grids. As a result, training data for MOR is best considered as graph-based data, necessitating methods that can handle such unstructured meshes. GNNs have arisen in recent years as powerful tools to learn on graphs, and have already celebrated remarkable successes in modelling complex physics problems \cite{Lam2023, Hernandez2022, Belbute2020}. These methods aim to enforce some geometric structure to a problem by embedding geometric priors into the model. Much work (reviewed in \cite{Zhang2018, Zhang2019, Liu2022}) has been dedicated to the the generalisation of pooling, unpooling and convolutional layers for graph-based data. GNNs are also not without limitations, with known issues for many graph convolutions such as reduced expressive power \cite{Danel2020} and oversmoothing \cite{Huang2020a, Li2018, Rusch2023}. In the MOR context, one seeks a map from a global parameter to a PDE solution, which can be represented as a graph. In GNN parlance, this represents a graph unpooling operation. However, currently there exist very few unpooling methods. In the major GNN libraries in Python, namely TF-GNN \cite{Ferludin2023}, PyG \cite{Fey2019} and DGL \cite{Wang2019}, there currently exists only one single graph unpooling layer, the $k$-nn interpolation layer, which is non-trainable and therefore clearly unsuitable for MOR.

\section{Graph Feedforward Network (GFN)}\label{sec:GFN}

As we discussed before, machine learning approaches for MOR usually exploit neural network layers which are not suitable for dealing with multifidelity data. Specifically, in the GNN context there is a lack of suitable unpooling layers for mapping from an input vector to an output graph, allowing for multi-resolution approaches. Given that feedforward networks have been heavily used in MOR, regardless of their potential inconsistencies, in this work we present GFNs, a generalisation of feedforward networks capable of leveraging data on any discretisation.

\begin{figure}[h]
    \centering
    \includegraphics[width=0.4\textwidth]{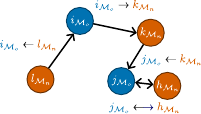}\hfill
    \caption{Illustration of the arrow notations used in this work.}
    \label{fig:notation}
\end{figure}

Prior to presenting the method, we fix some notation and conventions for ease of reading, also providing their graphical illustrations in \autoref{fig:notation}.

\begin{itemize}
    \item A mesh $\mathcal{M}$ is an ordered collection of nodes (co-ordinates), without any duplicates\footnote{We represent $\mathcal{M}$ as a set, whose ordering is arbitrary and not important here, provided it is fixed for a given set.}.
    \item $i_{\inmeshlg{}}$ denotes an index associated to a node in $\mathcal{M}$, defined by the ordering of the mesh $\mathcal{M}$. Thus, the position of a node of index $i_{\inmeshlg{}}$ is given as $\mathcal{M}\left[ i_{\inmeshlg{}} \right]$.
    \item $i_{\inmeshlg{o}} \leftarrow j_{\inmeshlg{n}}$ means that the node of index $i_{\inmeshlg{o}}$ in the mesh $\mathcal{M}_o$ is the nearest neighbour of the node of index $j_{\inmeshlg{n}}$ in the mesh $\mathcal{M}_n$, i.e., $i_{\inmeshlg{o}} = \operatorname{argmin}_{k_{\inmesh{o}}} \left|{\mathcal{M}_{n}} \left[ j_{\inmeshlg{n}} \right] - {\mathcal{M}_{o}}\left[ k_{\inmeshlg{o}} \right] \right|$.
    \item $i_{\inmeshlg{o}} \rightarrow j_{\inmeshlg{n}}$ means that the node of index $j_{\inmeshlg{n}}$ in the mesh $\mathcal{M}_n$ is the nearest neighbour\footnote{Note that we could write $i_{\inmeshlg{o}} \rightarrow j_{\inmeshlg{n}}$ as $j_{\inmeshlg{n}} \leftarrow i_{\inmeshlg{o}}$  and likewise $i_{\inmeshlg{o}} \leftarrow j_{\inmeshlg{n}}$ as $j_{\inmeshlg{n}} \rightarrow i_{\inmeshlg{o}}$. However, for readability we adopt the convention that nodes pertaining to $\mathcal{M}_o$ are given on the left side of the expression.} of the node of index $i_{\inmeshlg{o}}$ in the mesh $\mathcal{M}_o$, i.e., $j_{\inmeshlg{n}} = \operatorname{argmin}_{k_{\inmesh{n}}} \left|{\mathcal{M}_{n}} \left[ k_{\inmeshlg{n}} \right] - {\mathcal{M}_{o}}\left[ i_{\inmeshlg{o}} \right] \right|$.
    \item $i_{\inmeshlg{o}} \eitherarrow j_{\inmeshlg{n}}$ means that either $i_{\inmeshlg{o}} \leftarrow j_{\inmeshlg{n}}$ or $i_{\inmeshlg{o}} \rightarrow j_{\inmeshlg{n}}$.
    \item $i_{\inmeshlg{o}} \longleftrightarrow j_{\inmeshlg{n}}$ means that both $i_{\inmeshlg{o}} \leftarrow j_{\inmeshlg{n}}$ and $i_{\inmeshlg{o}} \rightarrow j_{\inmeshlg{n}}$.
\end{itemize}

\subsection{Feedforward networks for multifidelity data}
\begin{figure}[th]
    \centering
    \includegraphics[width=\textwidth]{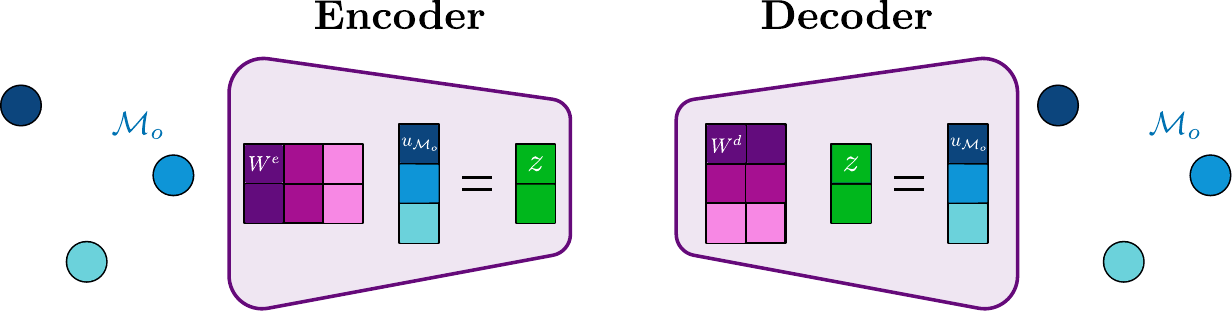}
    \caption{Single-layer graph feedforward autoencoder architecture. As shown via shading, the columns of the encoder weight matrix and the rows of the decoder weight matrix are associated to the individual nodes in the mesh $\mathcal{M}_o$.} 
    \label{fig:transforms_1}
\end{figure}

Numerous state-of-the-art ROMs such as GCA-ROM \cite{Pichi2023} and DL-ROM \cite{Fresca2021} employ autoencoder-based approaches for MOR.
In its simplest possible form, a standard feedforward network for an autoencoder task on an unstructured mesh $\mathcal{M}_o$ with $|\mathcal{M}_o|=N_o$ nodes consists of a single-layer encoder and single-layer decoder defined as
\bea
\operatorname{enc}(\bfmu_{\inmeshlg{o}})_i &=& \sigma \left( \sum_{j_{\inmesh{o}}=1}^{N_o} W^e_{ij_{\inmesh{o}}} u_{j_{\inmesh{o}}} + b^e_i \right), \qquad \forall i=1,\cdots,L,\\
\operatorname{dec}(\bfmz)_{i_{\inmesh{o}}} &=& \sum_{j=1}^L W^d_{i_{\inmesh{o}}j} z_j + b^d_{i_{\inmesh{o}}}, \qquad \forall i_{\inmeshlg{o}}=1,\cdots,N_o,
\eea
where $\bfmu_{\inmeshlg{o}} \in \mathbb{R}^{N_o}$ is an input sample on the mesh $\mathcal{M}_o$, $\bfmz \in \mathbb{R}^L$ is a latent vector, $\sigma$ is an activation function, $\bfmW^e \in \mathbb{R}^{L \times N_o}$ and $\bfmb^e \in \mathbb{R}^L$ are the encoder weights and biases, respectively, and $\bfmW^d \in \mathbb{R}^{N_o \times L}$ and $\bfmb^d \in \mathbb{R}^{N_o}$ are the decoder weights and biases, respectively. The reconstruction loss to optimise is simply given as $|| \operatorname{dec}\left(\operatorname{enc}(\bfmu_{\inmeshlg{o}}) \right) - \bfmu_{\inmeshlg{o}}||_2^2$, thus defining the machine learning task. We note that the encoder and decoder weights $W^e_{ji_{\inmesh{o}}}$, $W^d_{i_{\inmesh{o}}j}$ and the decoder bias $b^d_{i_{\inmesh{o}}}$ are all associated to a node $i_{\inmeshlg{o}}$. Consequently, we can regard these weights and biases as belonging to a node in the mesh $\mathcal{M}_o$, meaning that each has a spatial location given by the node's position. An illustration of this is shown for an autoencoder in \autoref{fig:transforms_1}.

\begin{figure}[t]
    \centering
    \includegraphics[width=\textwidth]{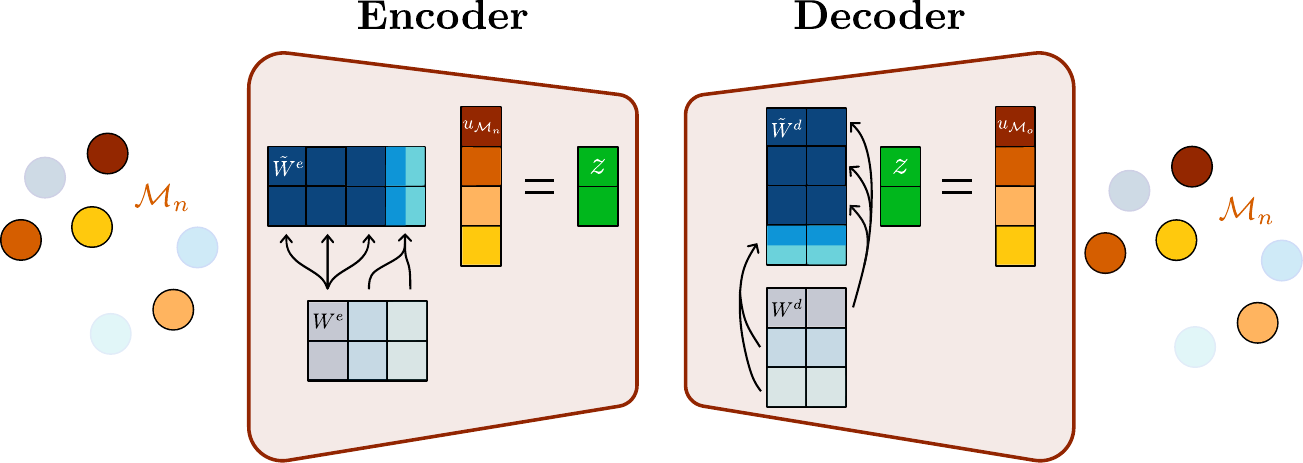}
    \caption{GFN approach to transform weights between meshes. Light blue shades denote the original mesh with its associated $\bfmW^e$ and $\bfmW^d$. The arrows depict how the new weight matrices (in darker blue shades) $\tilde{\bfmW}^e$ and $\tilde{\bfmW}^d$ can be created from the original weights and used for prediction on the new mesh.}
    \label{fig:transforms_2}
\end{figure}

Now, we suppose that we have a trained autoencoder defined above performing well on the unstructured mesh $\mathcal{M}_o$, and we wish to use this autoencoder on a new similar mesh $\mathcal{M}_n$ with $|\mathcal{M}_n|=N_n$ nodes. This is not possible for standard feedforward network, even if it is intuitively clear that one should be able to leverage the information from the previous model. As a result, we introduce the notion of a GFN and use the nearest neighbour information between the meshes $\mathcal{M}_o$ and $\mathcal{M}_n$ to generalise the autoencoder, with provable guarantees and desirable properties. Specifically, GFNs create a means of transforming weights and biases of the encoder and decoder $(\bfmW^e, ~\bfmb^e, ~\bfmW^d \text{ and } \bfmW^d)$ for the training mesh $\mathcal{M}_o$ to new weights and biases $( \tilde{\bfmW}^e, ~\tilde{\bfmb}^e, ~\tilde{\bfmW}^d \text{ and } \tilde{\bfmW}^d )$ for the new mesh $\mathcal{M}_n$, allowing for the evaluation of a feedforward network on multifidelity data. 
To achieve this transform from $\mathcal{M}_o$ to $\mathcal{M}_n$, which we denote coherently as $\mathcal{M}_o \rightarrow \mathcal{M}_n$, while retaining good performance, one can leverage the similarity of proximate nodes as illustrated in \autoref{fig:transforms_2}.

Mathematically, the GFN transforms for $\mathcal{M}_o \rightarrow \mathcal{M}_n$ are defined as
\begin{equation}
\begin{aligned}
\label{eq:gfn}
\tilde{W}^e_{ij_{\inmesh{n}}} &= \sum_{\forall k_{\inmesh{o}} \text{ s.t } k_{\inmesh{o}} \eitherarrow j_{\inmesh{n}}} \frac{W^e_{ik_{\inmesh{o}}}}{\lvert \{ h_{\inmeshlg{n}} \text{ s.t. } k_{\inmeshlg{o}} ~\eitherarrow~ h_{\inmeshlg{n}} \}\rvert}, \\
\tilde{b}^e_{i} &= b^e_{i}, \\
\tilde{W}^d_{i_{\inmesh{n}}j} &= \underset{\forall k_{\inmesh{o}} \text{ s.t } k_{\inmesh{o}} \eitherarrow i_{\inmesh{n}}}{\operatorname{mean}} {W}^d_{k_{\inmesh{o}}j}, \\
\tilde{b}^d_{i_{\inmesh{n}}} &=  \underset{\forall k_{\inmesh{o}} \text{ s.t } k_{\inmesh{o}} \eitherarrow i_{\inmesh{n}}}{\operatorname{mean}} {b}^d_{k_{\inmesh{o}}}.
\end{aligned}
\end{equation}
Note that trivially $\tilde{b}^e_i = b^e_i$ since the encoder bias is not associated with the graph, but rather with the latent vector. These transforms require nearest neighbour computations for all nodes in $\mathcal{M}_o$ and $\mathcal{M}_n$, leading to algorithm with time complexity of $O(N_o \log N_n + N_n \log N_o)$, using the $k$-d tree method for nearest neighbour computations \cite{Bentley1975}. The predictions on the new mesh $\mathcal{M}_n$ then read analogously to before as
\bea
\label{eq:enc_Mn}
\operatorname{enc}(\bfmu_{\inmeshlg{n}})_i &=& \sigma \left( \sum_{j_{\inmesh{n}}=1}^{N_n} \tilde{W}^e_{ij_{\inmesh{n}}} u_{j_{\inmesh{n}}} + b^e_i \right), \qquad \forall i=1,\cdots,L,\\
\label{eq:dec_Mn}
\operatorname{dec}(\bfmz)_{i_{\inmesh{n}}} &=& \sum_{j=1}^L \tilde{W}^d_{i_{\inmesh{n}}j} z_j + \tilde{b}^d_{i_{\inmesh{n}}}, \qquad \forall i_{\inmeshlg{n}}=1,\cdots,N_n,
\eea
where $\bfmu_{\inmeshlg{n}} \in \mathbb{R}^{N_n}$ is now a new input sample on the mesh $\mathcal{M}_n$ and the new weights and biases $\tilde{\bfmW}^e \in \mathbb{R}^{L \times N_n}$, $\tilde{\bfmb}^e \in \mathbb{R}^L$, $\tilde{\bfmW}^d \in \mathbb{R}^{N_n \times L}$ and $\tilde{\bfmb}^d \in \mathbb{R}^{N_n}$ are given in \autoref{eq:gfn}. Note that the transformed weights are now associated to each of the nodes on the new mesh $\mathcal{M}_n$.

To keep track of meshes and transforms, we denote the presented GFN approach of transforming weights from a mesh $\mathcal{M}_o$ to a new mesh $\mathcal{M}_n$ as 
\bea\label{eq:gfn_transform}
\tilde{\bfmW}^e, \tilde{\bfmW}^d, \tilde{\bfmb}^d = \operatorname{GFN}^{\mathcal{M}_o \to \mathcal{M}_n}(\bfmW^e, \bfmW^d, \bfmb^d).
\eea
Similarly, we also write $\operatorname{enc}^{\mathcal{M}_o \to \mathcal{M}_n}$ and $\operatorname{dec}^{\mathcal{M}_o \to \mathcal{M}_n}$ to track the applied transforms.

For autoencoders with multiple hidden layers, the GFN approach trivially generalises and only needs to be applied to the first encoder layer and the last decoder layer, since these are the only layers directly interacting with the changing meshes (for further details see \hyperref[app:mult_layers]{Appendix~\ref*{app:mult_layers}}).

This approach of transforming via GFN thus defines a means of allowing an autoencoder, which has been trained on single fidelity data using standard feedforward networks only, to be evaluated on multifidelity data.

\subsection{Algorithmic details}
\begin{figure}[h]
    \centering
    \includegraphics[width=.8\textwidth]{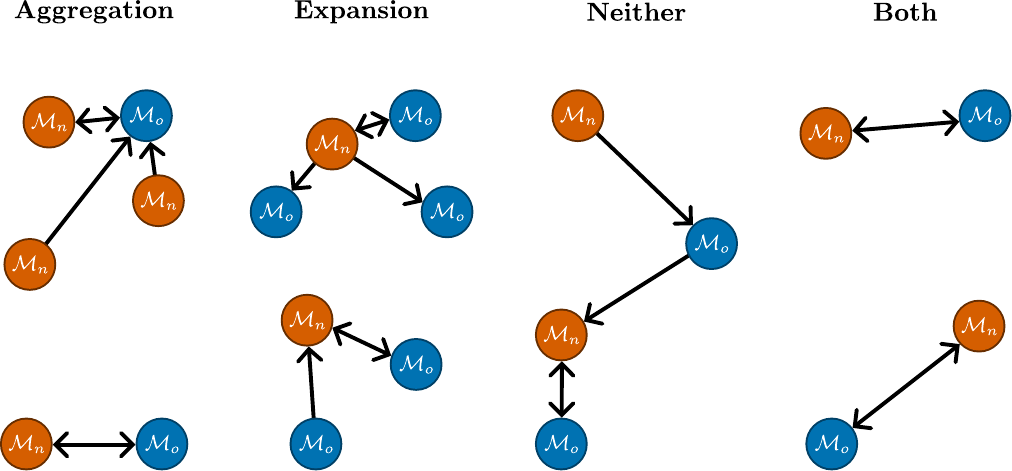}
    \caption{Examples of the different types of possible transforms.}
    \label{fig:example_agg_exp}
\end{figure}

In this section, we highlight two notable instances of the possible transforms $\mathcal{M}_o \to \mathcal{M}_n$ which lead to computational savings: the \textit{expansive} and the \textit{agglomerative} cases. Demonstrative examples of agglomerative, expansive, neither agglomerative nor expansive and both agglomerative and expansive cases are shown in \autoref{fig:example_agg_exp}.

\subsubsection{Expansive case}
The nearest neighbour of each node $i_{\inmeshlg{o}}$ in $\mathcal{M}_n$ has $i_{\inmeshlg{o}}$ as its nearest neighbour in $\mathcal{M}_o$, i.e., the nearest neighbour of the nearest neighbour of each node $i_{\inmeshlg{o}}$ is itself, following our notation:
    \bea
    \label{eq:exp}
    \forall i_{\inmeshlg{o}}, \quad i_{\inmeshlg{o}} \rightarrow j_{\inmeshlg{n}} \implies i_{\inmeshlg{o}} \leftarrow j_{\inmeshlg{n}}.
    \eea

In such case, the transform can be simplified (see the proof in \hyperref[subsec:exp_proof]{Appendix~\ref*{subsec:exp_proof}}) as:
    \begin{equation}
    \label{eq:gfn_exp}
    \begin{aligned}
    \tilde{W}^e_{ij_{\inmesh{n}}} &= \frac{1}{|\{h_{\inmeshlg{n}} \text{ s.t. } k_{\inmeshlg{o}} \leftarrow h_{\inmeshlg{n}} \}|} {W}^e_{ik_{\inmesh{o}}}, \qquad \text{ where}\quad k_{\inmeshlg{o}} \leftarrow j_{\inmeshlg{n}}, \\
    \tilde{W}^d_{i_{\inmesh{n}}j} &= {W}^d_{k_{\inmesh{o}}j}, \qquad \text{ where}\quad k_{\inmeshlg{o}} \leftarrow i_{\inmeshlg{n}}, \\
    \tilde{b}^d_{i_{\inmesh{n}}} &= {b}^d_{k_{\inmesh{o}}},  \qquad \text{ where}\quad k_{\inmeshlg{o}} \leftarrow i_{\inmeshlg{n}}.
    \end{aligned}
    \end{equation}
    We note that the absence of rightward arrows means we are only required to find nearest neighbours of nodes in the mesh $\mathcal{M}_n$ and not those in $\mathcal{M}_o$, reducing the algorithmic time complexity to $O(N_n \log N_o)$. The implementation of the expansive transformation with this complexity is shown in \hyperref[alg:gfn_expand]{Algorithm~\ref*{alg:gfn_expand}}.

\begin{algorithm}
\caption{Simplified GFN transform in the expansive case, as given in \autoref{eq:gfn_exp}.}
\label{alg:gfn_expand}
\begin{algorithmic}
\Function{Expand}{$\bfmW^e \in \mathbb{R}^{L \times N_o}, \bfmW^d \in \mathbb{R}^{N_o \times L}, \bfmb^d \in \mathbb{R}^{N_o}, \mathcal{M}_o, \mathcal{M}_n$}
\State{Initialise $\tilde{\bfmW}^e \in \mathbb{R}^{L \times N_n}, \tilde{\bfmW}^d \in \mathbb{R}^{N_n \times L}, \tilde{\bfmb}^d \in \mathbb{R}^{N_n}$}
\State $\operatorname{counts} \gets \boldsymbol{0} \in \mathbb{R}^{N_o}$
\For{$i_{\inmeshlg{n}}=1 \text{ to } N_n$}
\State $j_{\inmeshlg{o}} \gets \operatorname{nearest-neighbour}^{\inmeshlg{o}}(\mathcal{M}_n \left[ i_{\inmeshlg{n}} \right])$
\State $\operatorname{counts}[j_{\inmeshlg{o}}] \gets \operatorname{counts}[j_{\inmeshlg{o}}]+1$
\EndFor
\For{$i_{\inmeshlg{n}}=1 \text{ to } N_n$}
\State $j_{\inmeshlg{o}} \gets \operatorname{nearest-neighbour}^{\inmeshlg{o}}(\mathcal{M}_n \left[ i_{\inmeshlg{n}} \right])$
\State $\tilde{\bfmW}^e \left[:,~i_{\inmeshlg{n}}\right] \gets \bfmW^e \left[:,~j_{\inmeshlg{o}}\right] \mathbin{/} \operatorname{counts}[j_{\inmeshlg{o}}]$
\State $\tilde{\bfmW}^d \left[i_{\inmeshlg{n}},~:\right] \gets {\bfmW}^d \left[j_{\inmeshlg{o}},~:\right]$
\State $\tilde{\bfmb}^d \left[ i_{\inmeshlg{n}}\right] \gets{\bfmb}^d \left[ j_{\inmeshlg{o}}\right]$
\EndFor
\State \Return{$\tilde{\bfmW}^e, \tilde{\bfmW}^d, \tilde{\bfmb}^d$}
\EndFunction
\end{algorithmic}
\end{algorithm}

\subsubsection{Agglomerative case}
The nearest neighbour of each node $j_{\inmeshlg{n}}$ in $\mathcal{M}_o$ has $j_{\inmeshlg{n}}$ as its nearest neighbour in $\mathcal{M}_n$, i.e., the nearest neighbour of the nearest neighbour of each node $j_{\inmeshlg{n}}$ is itself, following our notation:
    \bea
    \label{eq:agg}
    \forall j_{\inmeshlg{n}}, \quad i_{\inmeshlg{o}} \leftarrow j_{\inmeshlg{n}} \implies i_{\inmeshlg{o}} \rightarrow j_{\inmeshlg{n}}.
    \eea
In such case, the transform can be simplified (see the proof in \hyperref[subsec:agg_proof]{Appendix~\ref*{subsec:agg_proof}}) as:
    \begin{equation}
    \label{eq:gfn_agg}
    \begin{aligned}
    \tilde{W}^e_{ij_{\inmesh{n}}} &=& \sum_{\forall k_{\inmesh{o}} \text{ s.t. }  k_{\inmesh{o}} \rightarrow j_{\inmesh{n}}} W^e_{ik_{\inmesh{o}}}, \\
    \tilde{W}^d_{i_{\inmesh{n}}j} &=& \underset{\forall k_{\inmesh{o}} \text{ s.t. }  k_{\inmesh{o}} \rightarrow i_{\inmesh{n}}}{\operatorname{mean}} W^d_{k_{\inmesh{o}}j}, \\
    \tilde{b}^d_{i_{\inmesh{n}}} &=&  \underset{\forall k_{\inmesh{o}} \text{ s.t. }  k_{\inmesh{o}} \rightarrow i_{\inmesh{n}}}{\operatorname{mean}} b^d_{k_{\inmesh{o}}}.
    \end{aligned}
    \end{equation}
    We note that the absence of leftward arrows means we are only required to find nearest neighbours of nodes in the mesh $\mathcal{M}_o$ and not those in $\mathcal{M}_n$, reducing the algorithmic time complexity to $O(N_o \log N_n)$. The implementation of the expansive transformation with this complexity is shown in \hyperref[alg:gfn_agglomerate]{Algorithm~\ref*{alg:gfn_agglomerate}}.

\begin{algorithm}
\caption{Simplified GFN transform in the agglomerative case, as given in \autoref{eq:gfn_agg}.}
\label{alg:gfn_agglomerate}
\begin{algorithmic}
\Function{Agglomerate}{$\bfmW^e \in \mathbb{R}^{L \times N_o}, \bfmW^d \in \mathbb{R}^{N_o \times L}, \bfmb^d \in \mathbb{R}^{N_o}, \mathcal{M}_o, \mathcal{M}_n$}
\State{Initialise $\tilde{\bfmW}^e \in \mathbb{R}^{L \times N_n}, \tilde{\bfmW}^d \in \mathbb{R}^{N_n \times L}, \tilde{\bfmb}^d \in \mathbb{R}^{N_n}$}
\State $\operatorname{counts} \gets \boldsymbol{0} \in \mathbb{R}^{N_n}$
\For{$i_{\inmeshlg{o}}=1 \text{ to } N_o$}
\State $j_{\inmeshlg{n}} \gets \operatorname{nearest-neighbour}^{\inmeshlg{n}}(\mathcal{M}_o \left[ i_{\inmeshlg{o}} \right])$
\State $\operatorname{counts}[j_{\inmeshlg{n}}] \gets \operatorname{counts}[j_{\inmeshlg{n}}]+1$
\State $\tilde{\bfmW}^e \left[:,~j_{\inmeshlg{n}}\right] \gets \bfmW^e \left[ :,~i_{\inmeshlg{o}}\right]$
\State $\tilde{\bfmW}^d \left[j_{\inmeshlg{n}},~:\right] \gets \left(\left( \operatorname{counts}[j_{\inmesh{n}}] - 1\right) \tilde{\boldsymbol{W}}^d \left[j_{\inmeshlg{n}},~:\right] + {\boldsymbol{W}}^d \left[i_{\inmeshlg{o}},~:\right]\right) \mathbin{/} \operatorname{counts}[j_{\inmeshlg{n}}]$
\State $\tilde{\bfmb}^d \left[ j_{\inmeshlg{n}}\right] \gets \left(\left( \operatorname{counts}[j_{\inmeshlg{n}}] - 1\right) \tilde{\boldsymbol{b}}^d \left[ j_{\inmeshlg{n}}\right] + {\boldsymbol{b}}^d \left[ i_{\inmeshlg{o}}\right]\right) \mathbin{/} \operatorname{counts}[j_{\inmeshlg{n}}]$
\EndFor
\State \Return{$\tilde{\bfmW}^e, \tilde{\bfmW}^d, \tilde{\bfmb}^d$}
\EndFunction
\end{algorithmic}
\end{algorithm}

\subsubsection{Computational advantages}
\label{subsubsec:comp_adv}
For hierarchical meshes, it follows that if $\mathcal{M}_o \subseteq \mathcal{M}_n$ then the transform $\mathcal{M}_o \rightarrow \mathcal{M}_n$ is expansive, whereas if $\mathcal{M}_n \subseteq \mathcal{M}_o$ then the transform $\mathcal{M}_o \rightarrow \mathcal{M}_n$ is agglomerative. Moreover, if the transform is both expansive and agglomerative then it must hold that $\forall i_{\inmeshlg{o}}$, $i_{\inmeshlg{o}} \longleftrightarrow j_{\inmeshlg{n}}$ and $\forall j_{\inmeshlg{n}}$, $i_{\inmeshlg{o}} \longleftrightarrow j_{\inmeshlg{n}}$, i.e.\ there is a one-to-one correspondence between nodes in the original mesh $\mathcal{M}_o$ and in the new mesh $\mathcal{M}_n$, and the GFN transforms reduce to nearest neighbour interpolation of the weights.

However, for general meshes a transform can be both, neither or one of expansive and agglomerative. Despite this, notably any general GFN transform can be expressed as a succession of an expansive update and an agglomerative update, (proof in \hyperref[app:gfn_decomp]{Appendix~\ref*{app:gfn_decomp}}) i.e.\ \autoref{eq:gfn} is equivalent to defining the new mesh
\bea
\label{eq:master_mesh}
\mathcal{M}_{} = \mathcal{M}_o \cup \{\mathcal{M}_n \left[i_{\inmeshlg{n}}\right] \text{ s.t. } j_{\inmeshlg{o}} \leftarrow i_{\inmeshlg{n}} \text{ but not } j_{\inmeshlg{o}} \rightarrow i_{\inmeshlg{n}} \},
\eea
and computing the transform in two steps, as
\begin{equation}
\begin{rcases}
\hat{W}^e_{ij_{\inmesh{}}} &= \frac{1}{|\{h_{\inmesh{}} \text{ s.t. } k_{\inmesh{o}} \leftarrow h_{\inmesh{}} \}|} {W}^e_{ik_{\inmesh{o}}}, \qquad \text{ where}\quad k_{\inmeshlg{o}} \leftarrow j_{\inmeshlg{}}, \\
\hat{W}^d_{i_{\inmesh{}}j} &= {W}^d_{k_{\inmesh{o}}j}, \qquad \text{ where}\quad k_{\inmeshlg{o}} \leftarrow i_{\inmeshlg{}}, \\
{{\hat{b}}^d}_{i_{\inmesh{}}} &= {b}^d_{k_{\inmesh{o}}},  \qquad \text{ where}\quad k_{\inmeshlg{o}} \leftarrow i_{\inmeshlg{}},
\end{rcases} \quad \text{Expansion step}
\end{equation}
followed by
\begin{equation}
\begin{rcases}
\tilde{W}^e_{ij_{\inmesh{n}}} &= \sum_{\forall k_{\inmesh{}} \text{ s.t. }  k_{\inmesh{}} \rightarrow j_{\inmesh{n}}} {\hat{W}}^e_{ik_{\inmesh{}}}, \\
\tilde{W}^d_{i_{\inmesh{n}}j} &= \underset{\forall k_{\inmesh{}} \text{ s.t. }  k_{\inmesh{}} \rightarrow i_{\inmesh{n}}}{\operatorname{mean}} {\hat{W}}^d_{k_{\inmesh{}}j}, \\
\tilde{b}^d_{i_{\inmesh{n}}} &=  \underset{\forall k_{\inmesh{}} \text{ s.t. }  k_{\inmesh{}} \rightarrow i_{\inmesh{n}}}{\operatorname{mean}} {\hat{b}}^d_{k_{\inmesh{}}}.
\end{rcases} \quad \text{Agglomerative step}
\end{equation}

\begin{figure}[ht]
    \centering
    \includegraphics[width=.8\textwidth]{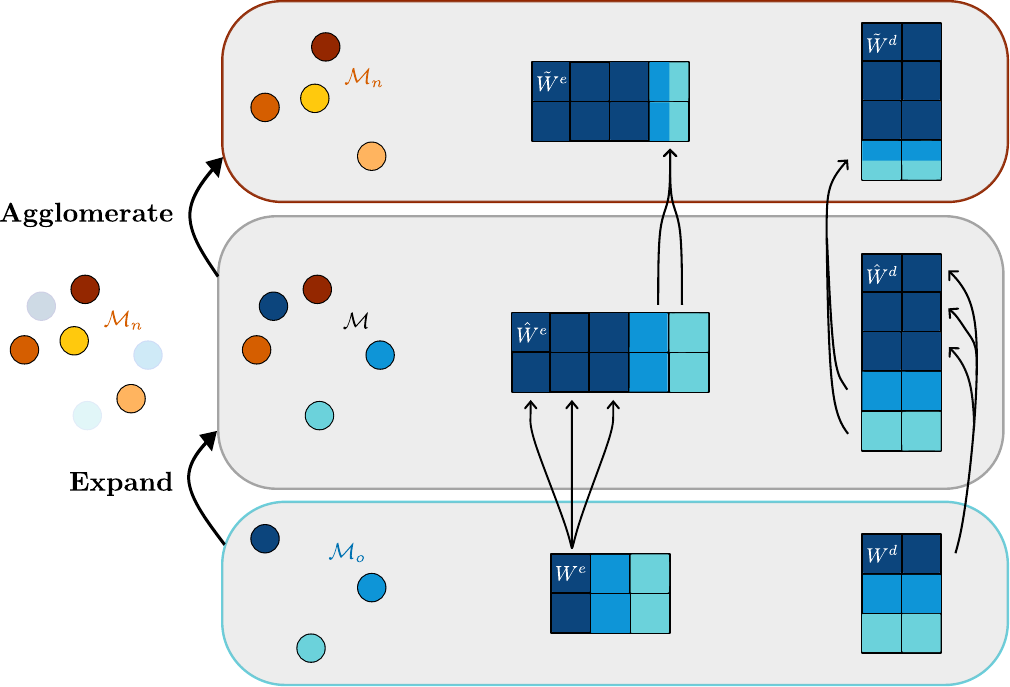}
    \caption{The GFN transform can be thought about as applying an expansive update followed by an agglomerative update.}
    \label{fig:gfn_exp_agg}
\end{figure}

The general GFN algorithm can thus be easily expressed as in \hyperref[alg:gfn_general]{Algorithm~\ref*{alg:gfn_general}}, with a demonstrative illustration of this interpretation in \autoref{fig:gfn_exp_agg}.

\begin{algorithm}
\caption{General GFN transform leveraging an expansion step followed by an agglomerative update, as given in \autoref{eq:gfn}.}
\label{alg:gfn_general}
\begin{algorithmic}
\Function{GFN}{$\bfmW^e \in \mathbb{R}^{L \times N_o}, \bfmW^d \in \mathbb{R}^{N_n \times L}, \bfmb^d \in \mathbb{R}^{N_o}, \mathcal{M}_o, \mathcal{M}_n$}
\State $\mathcal{M} = \mathcal{M}_o \cup \{\mathcal{M}_n \left[i_{\inmeshlg{n}}\right] \text{ s.t. } j_{\inmeshlg{o}} \leftarrow i_{\inmeshlg{n}} \text{ but not } j_{\inmeshlg{o}} \rightarrow i_{\inmeshlg{n}} \}$
\State $\hat{\bfmW}^e, \hat{\bfmW}^d, \hat{\bfmb}^d \gets \Call{Expand}{\bfmW^e, \bfmW^d, \bfmb^d, \mathcal{M}_o, \mathcal{M}}$
\State $\tilde{\bfmW}^e, \tilde{\bfmW}^d, \tilde{\bfmb}^d \gets \Call{Agglomerate}{\hat{\bfmW}^e, \hat{\bfmW}^d, \hat{\bfmb}^d, \mathcal{M}, \mathcal{M}_n}$
\State \Return{$\tilde{\bfmW}^e, \tilde{\bfmW}^d, \tilde{\bfmb}^d$}
\EndFunction
\end{algorithmic}
\end{algorithm}

\section{Graph Feedforward Network Reduced Order Model (GFN-ROM)}\label{sec:GFN-ROM}
\begin{figure}[t]
    \centering
\includegraphics[width=0.9\textwidth]{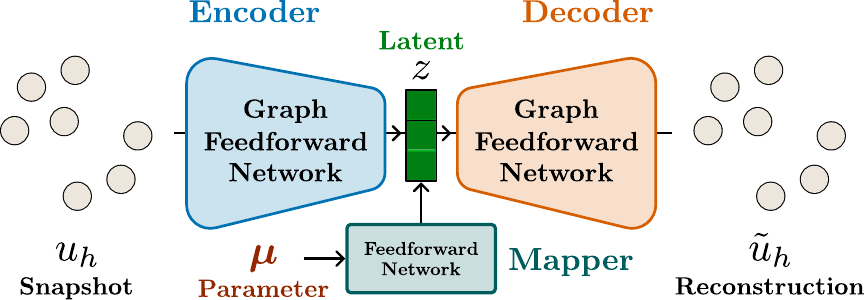}
    \caption{GFN-ROM architecture, consisting of an autoencoder and a mapper. During inference, predictions can be computed by decoding the mapper output.}
    \label{fig:gfn_rom}
\end{figure}

Having defined GFNs as a means of evaluating trained autoencoders on multifidelity data, we leverage them to propose a reduced order model capable of being evaluated on arbitrary meshes.
As depicted in \autoref{fig:gfn_rom}, the proposed GFN-ROM architecture consists of two components: an autoencoder and a mapper. The autoencoder is based on GFNs, while the mapper can be any neural network model which maps the fixed-size vector parameters directly to the fixed-size latent representation. Notice that, within this parametric context, the mapper does not need to be defined using GFNs since it does not deal with graph-structured data.

GFN-ROM is inspired by recent state-of-the-art autoencoder-based models such as DL-ROM \cite{Fresca2021}, which uses standard convolutions, and GCA-ROM \cite{Pichi2023}, which uses graph convolutions. In fact, GFN-ROM can be regarded as an extension of such methods. Like both models, we train the GFN-ROM model based on a weighted sum of two losses: the reconstruction loss and the mapper loss. 
Given a GFN-ROM model with weights defined on a mesh $\mathcal{M}_o$, the per-example loss functions are defined for the $t$th training sample $\left( \boldsymbol{\mu}_t, ~\bfmu_{\inmeshlgtwo{n}{t}}(\boldsymbol{\mu}_t) \right)$, defined on the $t$-th mesh $\mathcal{M}_n^t$ and for the $t$-th parameter value $\boldsymbol{\mu}_t$, as
\bea
\mathcal{L}_{\text{recon}}[\boldsymbol{\mu}_t,\bfmu_{\inmeshlgtwo{n}{t}}(\boldsymbol{\mu}_t)] &=& \frac{1}{|\mathcal{M}_n^t|} || \operatorname{dec}^{\mathcal{M}_o \to \mathcal{M}_n^t}\left(\operatorname{enc}^{\mathcal{M}_o \to \mathcal{M}_n^t}(\bfmu_{\inmeshlgtwo{n}{t}}(\boldsymbol{\mu}_t)) \right) - \bfmu_{\inmeshlgtwo{n}{t}}(\boldsymbol{\mu}_t) ||_2^2, \\
\mathcal{L}_{\text{map}}[\boldsymbol{\mu}_t, \bfmu_{\inmeshlgtwo{n}{t}}(\boldsymbol{\mu}_t)] &=& \frac{1}{L} || \operatorname{enc}^{\mathcal{M}_o \to \mathcal{M}_n^t}(\bfmu_{\inmeshlgtwo{n}{t}}(\boldsymbol{\mu}_t)) - \operatorname{map}(\boldsymbol{\mu}_t)||_2^2.
\eea
The overall loss is given as a weighted mean of these losses over a set $\{ \left( \boldsymbol{\mu}_t, ~ \bfmu_{\inmeshlgtwo{n}{t}}(\boldsymbol{\mu}_t) \right) \}_{t=1}^{T}$ of $T$ training examples 
\bea
\mathcal{L} = \frac{1}{T} \sum_{t=1}^{T} \frac{|\mathcal{M}_n^t|}{\sum_{s=1}^{T} |\mathcal{M}_n^s|} \left( \mathcal{L}_{\text{recon}}[\boldsymbol{\mu}_t,\bfmu_{\inmeshlgtwo{n}{t}}(\boldsymbol{\mu}_t)] + \omega \mathcal{L}_{\text{map}}[\boldsymbol{\mu}_t,\bfmu_{\inmeshlgtwo{n}{t}}(\boldsymbol{\mu}_t)] \right),
\eea
where we have introduced the mapper weight hyperparameter $\omega$.
We note that at inference, the encoder is no longer required and, for a parameter $\boldsymbol{\mu}_t$ and mesh $\mathcal{M}_n^t$, a GFN-ROM model with weights defined on a mesh $\mathcal{M}_o$ predicts $\operatorname{dec}^{\mathcal{M}_o \to \mathcal{M}_n^t}(\operatorname{map}(\boldsymbol{\mu}_t))$.

\subsection{Main features of the architecture}
The proposed GFN-ROM architecture benefits from a number of desirable properties that we discuss in the following sections, which make it an attractive tool for applications with multifidelity data.

\subsubsection{Self-consistency} The GFN transforms are self-consistent approaches, meaning that a series of two transforms from $\mathcal{M}_o$ to $\mathcal{M}_n$ and then back to $\mathcal{M}_o$ recovers exactly the original weights $\bfmW^e, \bfmW^d, \bfmb^d$ provided that the transform $\mathcal{M}_o$ to $\mathcal{M}_n$ is expansive, i.e.\ $$\bfmW^e, \bfmW^d, \bfmb^d = \operatorname{GFN}^{\mathcal{M}_n \to \mathcal{M}_o}( \operatorname{GFN}^{\mathcal{M}_o \to \mathcal{M}_n}(\bfmW^e, \bfmW^d, \bfmb^d) ).$$
The self-consistency, a proof of which is given in \hyperref[app:self_consistency]{Appendix~\ref*{app:self_consistency}}, is a marked advantage in comparison to any other interpolation methods such as Kriging, linear regression, spline, and $k$NN (with the exception of the case $k=1$). This means that no information is lost during the expansive transform and the original weights can be directly recovered, unlike for other interpolation-based approaches.

\subsubsection{Bound on super- and sub-resolution error of GFN-ROM}
GFN-ROM is a resolution-invariant architecture, meaning one is often interested in evaluating the model on a finer discretisation (super-resolution) or coarser discretisation (sub-resolution) compared to a reference resolution. A particularly attractive feature of the architecture is that one can bound the performance on any testing fidelity given certain information. Specifically, given a bound $\delta$ on the difference between the features of nearest neighbours, and a bound $\tau$ on the performance of GFN-ROM on the known mesh $\mathcal{M}_o$, i.e.\footnote{We use the notation $u(x_{i_{\inmesh{o}}})$, equivalent to $u_{i_{\inmesh{o}}}$, in order to emphasise the dependence of $u$ on spatial position.}
\bea
\forall i_{\inmeshlg{o}}, \quad i_{\inmeshlg{o}} \eitherarrow j_{\inmeshlg{n}} \implies \left| u(x_{i_{\inmesh{o}}} ) - u(x_{j_{\inmesh{n}}}) \right| \leq \delta, \label{eqn:mesh_bound} \\
\forall i_{\inmeshlg{o}},\quad \left| u(x_{i_{\inmesh{o}}}) - \operatorname{dec}^{\mathcal{M}_o \to \mathcal{M}_o}(\operatorname{map}(\boldsymbol{\mu}))_{i_{\inmesh{o}}} \right| \leq \tau, \label{eqn:autoencoder_bound_1}
\eea
one can show that the performance of GFN-ROM on the new mesh $\mathcal{M}_n$ can be bounded as:
\bea
\forall i_{\inmeshlg{n}},\quad \left| u(x_{i_{\inmesh{n}}}) - \operatorname{dec}^{\mathcal{M}_o \to \mathcal{M}_n}(\operatorname{map}(\boldsymbol{\mu}))_{i_{\inmesh{n}}} \right| \leq \tau + \delta.
\eea
A proof is given in \hyperref[subsec:error_gfnrom]{Appendix~\ref*{subsec:error_gfnrom}}. We note that there is no dependence at all on the weights and biases, only on the mesh similarities and the GFN-ROM performance on the reference resolution.

\subsubsection{Bound on super- and sub-resolution error of the GFN-ROM mapper}
Given the same bound $\delta$ as in \autoref{eqn:mesh_bound}, and a bound $\alpha$ on the performance of the mapper on the known mesh $\mathcal{M}_o$, i.e.\
\bea
\forall i,\quad \left| \operatorname{map}(\boldsymbol{\mu})_i - \operatorname{enc}^{\mathcal{M}_o \to \mathcal{M}_o}(\bfmu_{\inmeshlg{o}})_i \right| \leq \alpha, \label{eqn:autoencoder_bound_2}
\eea
one can show that the performance of the mapper on the new mesh $\mathcal{M}_n$ can be bounded as:
\bea
\forall i,\quad \left| \operatorname{map}(\boldsymbol{\mu})_i - \operatorname{enc}^{\mathcal{M}_o \to \mathcal{M}_n}(\bfmu_{\inmeshlg{n}})_i \right| \leq \alpha + \delta C^{P+1} ||\bfmW^{e}||_\infty.
\eea
A proof is given in \hyperref[subsec:error_mapper]{Appendix~\ref*{subsec:error_mapper}}.

\subsubsection{Bound on super- and sub-resolution error of the GFN-ROM autoencoder}
Given the same bound $\delta$ as in \autoref{eqn:mesh_bound}, and a bound $\beta$ on the performance of the autoencoder on the known mesh $\mathcal{M}_o$, i.e.\
\bea
\forall i_{\inmeshlg{o}}, \left| u(x_{i_{\inmesh{o}}}) - \operatorname{dec}^{\mathcal{M}_o \to \mathcal{M}_o}(\operatorname{enc}^{\mathcal{M}_o \to \mathcal{M}_o}(\bfmu_{\inmeshlg{o}}))_{i_{\inmesh{o}}} \right| \leq \beta, \label{eqn:autoencoder_bound_3}
\eea
one can show that the performance of the autoencoder on the new mesh $\mathcal{M}_n$ can be bounded as:
\begin{equation}
\begin{aligned}
\forall i_{\inmeshlg{n}},\quad \left| u(x_{i_{\inmesh{n}}}) - \operatorname{dec}^{\mathcal{M}_o \to \mathcal{M}_n}(\operatorname{enc}^{\mathcal{M}_o \to \mathcal{M}_n}(\bfmu_{\inmeshlg{n}}))_{i_{\inmesh{n}}} \right| \leq \beta + \delta +  \delta C^{Q+1} ||\bfmW^{d}||_\infty ||\bfmW^{e}||_\infty.
\end{aligned}
\end{equation}
A proof is given in \hyperref[subsec:error_autoenc]{Appendix~\ref*{subsec:error_autoenc}}. The dependence on the norm of the weights in the bound suggests that regularisation is important to allow for good generalisation for the autoencoder.

\subsection{Link between GFN-ROM and novel architectures}
Following the discussion in \hyperref[subsec:nos]{Section~\ref*{subsec:nos}} and \hyperref[subsec:gnns]{Section~\ref*{subsec:gnns}}, we discuss here how GFN-ROM is connected with state-of-the-art resolution invariant approaches.

\subsubsection{GNNs and GFN-ROM}
Concerning graph neural networks, we can interpret a GFN encoder layer and a GFN decoder layer as two different types of GNN: namely, a graph pooling layer and a graph unpooling layer, respectively. We can consider the input to the autoencoder, which is some PDE solution $u(x, \boldsymbol{\mu})$ evaluated on some discretized domain, as a fully connected graph. The GFN encoder layer is then tasked with computing a vector representation summarising the input graph, i.e.\ a graph pooling task. The GFN decoder layer has precisely the opposite task of reconstructing a graph given a summarised representation i.e.\ graph unpooling. At inference, the overall GFN-ROM architecture can thus be thought of as a graph unpooling model, where the PDE solution is identified by the parameters $\boldsymbol{\mu}$.

\subsubsection{NOs and GFN-ROM}
A parameterised PDE with a parameter $a \in \mathcal{A}$ and solution $u \in \mathcal{U}$, where $\mathcal{A}$ and $\mathcal{U}$ are Banach spaces, defines a mapping $\mathcal{G}^\dag: \mathcal{A} \to \mathcal{U}$, i.e.\ $u(x)=\mathcal{G}^\dag[a](x)$. A neural operator aims to directly approximate this mapping $\mathcal{G}^\dag$ by means of an operator $\mathcal{G}_\theta \approx \mathcal{G}^\dag$ with trainable parameters $\theta$.
A single-layer neural operator can be written as
\nbea
G_\theta = \mathcal{Q} \circ \sigma(W+\mathcal{K} + b) \circ \mathcal{P},
\neea
where $\mathcal{P}$, $\mathcal{Q}$ are the local lifting and projection mappings, $W$ is a local linear operator (matrix), $\mathcal{K}$ is an integral kernel operator, $b$ is a bias function and $\sigma$ is a pointwise activation function. We refer to Kovachki \emph{et al.} for further details \cite{Kovachki2021}. Considering a simplified neural operator with $W=0$ and no lifting or projection mappings, we obtain
\bea
\label{eq:no}
G_\theta[v](y) = \sigma \left( \int_D \kappa(x,y) v(x)\, \mathrm{d}x + b(y) \right),
\eea
where $\kappa$ is a kernel function and $D$ is the domain of integration. Usually, one does not have complete access to the function $v(x)$, but only to its evaluations at certain points, meaning this integral has to be numerically computed. Supposing we are given function evaluations on a mesh $\mathcal{M}$, the numerical approximation of the integral reads as
\bea
G_\theta[v](y) \approx \sigma\left(\sum_{i_{\inmesh{}}} h_{i_{\inmesh{}}} \kappa(x_{i_{\inmesh{}}},y) v(x_{i_{\inmesh{}}}) + b(y) \right),
\eea
where $h_{i_{\inmesh{}}}$ is a factor depending on the integration method.

GFN-ROM is closely related to NOs. In fact, in the expansive case, GFN-ROM can be recast as a particular instance of NOs, since we can directly obtain both the encoder and decoder. To make the link with GFN-ROM, we start with the mesh $\mathcal{M}_o$, the associated weights and biases ($\bfmW^e$, $\bfmb^e$, $\bfmW^d$, $\bfmb^d$), and choose a particular kernel $\kappa$.

Starting from \autoref{eq:no}, discarding the $y$-dependence and taking $D=\Omega$, $b=b^e_i$, and $\kappa(x) = W^e_{ik_{\inmesh{o}}}/|B_{k_{\inmesh{o}}}|$ where $k_{\inmeshlg{o}} \leftarrow x$, and $B_{k_{\inmesh{o}}} = \{x \in \Omega\ |\ k_{\inmeshlg{o}} \leftarrow x \}$, one can write\footnote{We abuse notation to write $k_{\inmeshlg{o}} \leftarrow x$ since $x$ is a position and not an index, but the meaning is analogous to before, i.e.\ $k_{\inmeshlg{o}} = \operatorname{argmin}_{h_{\inmesh{o}}} \left|{\mathcal{M}_{o}}\left[ h_{\inmeshlg{o}} \right]  - x \right|$.}
\begin{equation}
\begin{aligned}
G_\theta[u] &= \sigma \left( \int_D \kappa(x) u(x) \, \mathrm{d}x + b \right),\\
&= \sigma \left( \int_\Omega \frac{W^e_{ik_{\inmesh{o}}}}{|B_{k_{\inmesh{o}}}|} u(x) \, \mathrm{d}x  + b^e_i \right), \qquad k_{\inmeshlg{o}} \leftarrow x,\\
&= \sigma \left( \sum_{k_{\inmesh{o}}=1}^{N_o} W^e_{ik_{\inmesh{o}}} \frac{1}{|B_{k_{\inmesh{o}}}|}\int_{B_{k_{\inmesh{o}}}} u(x) \, \mathrm{d}x + b^e_i \right).
\end{aligned}
\end{equation}
Approximating the integral with observations of $u(x)$ on a mesh $\mathcal{M}_n$, with $k_{\inmeshlg{o}} \leftarrow j_{\inmeshlg{n}}$, gives the formula for the GFN-ROM encoder in the expansive case
\begin{equation}
\begin{aligned}
G_\theta[u] &= \sigma \left( \sum_{k_{\inmesh{o}}=1}^{N_o} W^e_{ik_{\inmesh{o}}} ~\underset{\forall j_{\inmesh{n}} \text{ s.t } k_{\inmesh{o}} \leftarrow j_{\inmesh{n}}}{\operatorname{mean}} u_{j_{\inmesh{n}}} + b^e_i \right), \\
&= \sigma \left( \sum_{j_{\inmesh{n}}=1}^{N_n} \frac{W^e_{ik_{\inmesh{o}}}}{\lvert \{ h_{\inmeshlg{n}} \text{ s.t. } k_{\inmeshlg{o}} \leftarrow h_{\inmeshlg{n}} \}\rvert}  u_{j_{\inmesh{n}}} + b^e_i \right),  \qquad \text{ where}\quad k_{\inmeshlg{o}} \leftarrow j_{\inmeshlg{n}},
\\
&= \operatorname{enc}^{\mathcal{M}_o \to \mathcal{M}_n}(\bfmu_{\inmeshlg{n}})_i.
\end{aligned}
\end{equation}
Similarly, the decoder can be obtained from \autoref{eq:no} without any $x$-dependence. 
Taking $\tilde{\kappa} = |D| \kappa$, setting $\sigma$ to the identity function, $\tilde{\kappa}(y)=\bfmW^d_{k_{\inmesh{o}},:}$ where $k_{\inmeshlg{o}} \leftarrow y$, and $b(y)=b^d_{k_{\inmesh{o}}}$ where $k_{\inmeshlg{o}} \leftarrow y$, one obtains
\begin{equation}
\begin{aligned}
G_\theta[\bfmz](y) &= \sigma \left( \int_D \kappa(y) \bfmz \, \mathrm{d}x + b(y) \right), \\
&= \sigma \left( \tilde{\kappa}(y) \bfmz + b(y) \right),\\
&= \sum_{j=1}^L {W}^d_{k_{\inmesh{o}}j} z_j + {b}^d_{k_{\inmesh{o}}}, \qquad k_{\inmeshlg{o}} \leftarrow y.
\end{aligned}
\end{equation}
Considering $y=\mathcal{M}_n[i_{\inmeshlg{n}}]$, we recover the GFN-ROM decoder in the expansive case as
\begin{equation}
\begin{aligned}
G_\theta[\bfmz](y) &= \sum_{j=1}^L {W}^d_{k_{\inmesh{o}}j} z_j + {b}^d_{k_{\inmesh{o}}},\\
&= \operatorname{dec}^{\mathcal{M}_o \to \mathcal{M}_n}(\bfmz)_{i_{\inmesh{n}}}.
\end{aligned}
\end{equation}
Like the decoder, the GFN-ROM model at inference in the expansive case can be thought of as a neural operator, but with the mapper acting as the lifting function. Thus, in such settings, GFN-ROM encoder and decoder are very closely linked to resolution-invariant neural operators with a choice of piecewise constant kernel function.

\subsubsection{ROM-related properties}
In addition to the above proofs, the GFN-ROM architecture benefits from a number of typical properties important for reduced order models. Namely, it is purely data-driven and non-intrusive, not requiring any knowledge of the generation process for the training data. It is furthermore a nonlinear method, leveraging an autoencoder for nonlinear compression. The method is suitable for unstructured meshes, which is vital in order to deal with general PDEs, where the presence of complex domains typically lead to unstructured meshes. Most importantly, the method is a completely multifidelity approach and is capable of inference on arbitrary meshes.

\section{Adaptive Multifidelity Training}\label{sec:adaptive_training}
The GFN method proposes a means of transferring weights from one graph $\mathcal{M}_o$ to a new graph $\mathcal{M}_n$ as $\tilde{\bfmW}^e, \tilde{\bfmW}^d, \tilde{\bfmb}^d = \operatorname{GFN}^{\mathcal{M}_o \to \mathcal{M}_n}(\bfmW^e, \bfmW^d, \bfmb^d)$, allowing one to evaluate a feedforward-based approach on new meshes. Whilst this is already a major advantage, it does not prescribe how to effectively train in a multifidelity setting.

\begin{figure}[ht]
    \centering
    \includegraphics[width=0.7\textwidth]{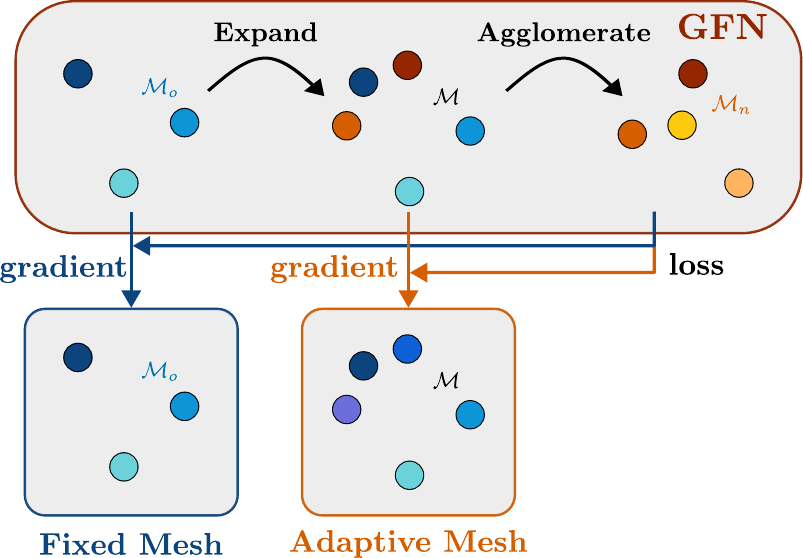}
    \caption{Training with  fixed or adaptive mesh methods when given a starting mesh $\mathcal{M}_o$ and new mesh $\mathcal{M}_n$. The fixed mesh approach updates the original mesh weights, whilst the adaptive mesh approach updates weights on a new mesh $\mathcal{M}$, created from information from both meshes.}
    \label{fig:adaptive_method}
\end{figure}

In this section, we discuss means of training in multifidelity applications. We aim at learning from a series of $T$ meshes $\mathcal{M}_n^1,\cdots,\mathcal{M}_n^T$ by training the GFN-ROM model with weights and biases corresponding to a fixed or an adaptive mesh $\mathcal{M}_o$. We describe below the two approaches, the difference between which is illustrated in \autoref{fig:adaptive_method}.

\subsection{Fixed mesh} The GFN transforms are entirely differentiable, meaning that for weights and biases $\bfmW^e$, $\bfmW^d$ and $\bfmb^d$ associated on a fixed mesh $\mathcal{M}_o$, it is possible to compute gradients $\nabla_{\boldsymbol{W}^e} \mathcal{L}$, $\nabla_{\boldsymbol{W}^d} \mathcal{L}$ and $\nabla_{\boldsymbol{b}^d} \mathcal{L}$ when training on new meshes. As a result, it is in fact possible to use the standard gradient-based optimisation methods to learn the optimal weights $\bfmW^e$, $\bfmW^d$ and $\bfmb^d$, which are always associated to the mesh $\mathcal{M}_o$. If one considers the mesh $\mathcal{M}_o$ as being a suitable sampling of the domain, this approach therefore suffices as a means of training the model in a multifidelity setting. However, the choice of mesh $\mathcal{M}_o$ is extremely important for the performance of the model. The fixed nature of $\mathcal{M}_o$ is also undesirable since, like interpolation approaches onto fixed grids, it does not allow the model to adapt the mesh $\mathcal{M}_o$ to better match the training data it sees.

\subsection{Adaptive mesh}
GFN-ROM differs from interpolation onto a fixed grid since it is possible to transform weights and biases from a coarse mesh into weights and biases for a finer mesh. Thus, one can always overwrite the old weights and biases corresponding to the old mesh with new weights and biases of the new mesh, and then optimise the new weights. As a result, the approach leads to a number of advantages in allowing the architecture to adapt its number of parameters during training. That is, instead of having a single fixed mesh $\mathcal{M}_o$ associated to the fixed weights and biases for the model, we replace it by a master mesh $\mathcal{M}^t$ which can change during each training sample $t=1,\cdots,T$. Importantly, the sequence of master meshes during training must be hierarchically increasing i.e. $\mathcal{M}^0 \subseteq \cdots \subseteq \mathcal{M}^T$, since we do not wish to lose information at any point.

A simple scheme for adaptive training could be setting $\mathcal{M}^t = \mathcal{M}^{t-1} \cup \mathcal{M}_n^t$. However, this approach is suboptimal since it necessarily adds every node seen during training to the master mesh. This means that the master mesh can become prohibitively large, and result in a worse time complexity of $O((N_o+N_n)\log N_o + (N_o+N_n)\log N_n)$ for the transform $\mathcal{M}_o \to \mathcal{M}_n$.

Instead, a better approach is possible without incurring the computational issues of the previous algorithm. As discussed in \autoref{subsubsec:comp_adv} and summarised in \hyperref[alg:gfn_general]{Algorithm~\ref*{alg:gfn_general}}, the GFN transforms can be decomposed into an expansive and agglomerative part via the choice of a suitable intermediate mesh given in \autoref{eq:master_mesh}. This intermediate mesh is already computed during the GFN transforms, meaning it does not represent an extra operation. Furthermore, the intermediate mesh is larger than the original mesh. Therefore, it is possible to use such mesh to develop a computationally efficient adaptive approach i.e., $\mathcal{M}^t = \mathcal{M}^{t-1} \cup  \{\mathcal{M}_n^t \left[i_{\inmeshlgtwo{n}{t}}\right] \text{ s.t. } j_{\inmeshlgtwo{}{t-1}} \leftarrow i_{\inmeshlgtwo{n}{t}} \text{ but not } j_{\inmeshlgtwo{}{t-1}} \rightarrow i_{\inmeshlgtwo{n}{t}} \}$. 

For a transform $\mathcal{M}_o \to \mathcal{M}_n$, we retain the time complexity of the non-adaptive algorithm. Furthermore, this approach only adds nodes to the master mesh when the master mesh undersamples\footnote{Undersampling here refers purely to nearest-neighbour information i.e., areas where a node in the master mesh has multiple nearest neighbours are considered undersampled.}. We note that because the weights and biases are changing shape along with the master mesh during training for adaptive methods, one cannot use Adam, RMSProp or any momentum-based optimisation methods for updates\footnote{In practice, we often know all of the training fidelities we shall see during training. In such cases, the final mesh $\mathcal{M}^T$ can be precomputed and the model trained as for a fixed mesh.}.

The training loop for GFN-ROM is given in \hyperref[alg:gfn_train]{Algorithm~\ref*{alg:gfn_train}}, illustrating the difference between the adaptive and fixed mesh modes.

\begin{algorithm}
\caption{Training algorithm for GFN-ROM, showing both adaptive and fixed mesh approaches.}
\label{alg:gfn_train}
\begin{algorithmic}
\Require $T$ training samples of parameter and solution pairs i.e. $\{ \left( \boldsymbol{\mu}_t,~ \bfmu_{\inmeshlgtwo{n}{t}}(\boldsymbol{\mu}_t) \right) \}_{t=1}^{T}$
\State Initialise a mesh $\mathcal{M}_o$ with associated GFN weights ${\bfmW}^e, {\bfmW}^d, {\bfmb}^d$
\State Initialise $\theta$: all parameters needed for GFN-ROM excluding ${\bfmW}^e, {\bfmW}^d, {\bfmb}^d$
\For{$t\gets 1, \dots, T$}
\State $\mathcal{M} = \mathcal{M}_o \cup \{\mathcal{M}^t_n \left[i_{\inmeshlgtwo{n}{t}}\right] \text{ s.t. } j_{\inmeshlg{o}} \leftarrow i_{\inmeshlgtwo{n}{t}} \text{ but not } j_{\inmeshlg{o}} \rightarrow i_{\inmeshlgtwo{n}{t}} \}$
\If{adaptive} \Comment{adaptive mode}
\State ${\bfmW}^e, {\bfmW}^d, {\bfmb}^d \gets \Call{Expand}{\bfmW^e, \bfmW^d, \bfmb^d, \mathcal{M}_o, \mathcal{M}}$
\State $\mathcal{M}_o \gets \mathcal{M}$
\State $\tilde{\bfmW}^e, \tilde{\bfmW}^d, \tilde{\bfmb}^d \gets \Call{Agglomerate}{{\bfmW}^e, {\bfmW}^d, {\bfmb}^d, \mathcal{M}_o, \mathcal{M}_n}$
\Else \Comment{fixed mode}
\State $\hat{\bfmW}^e, \hat{\bfmW}^d, \hat{\bfmb}^d \gets \Call{Expand}{\bfmW^e, \bfmW^d, \bfmb^d, \mathcal{M}_o, \mathcal{M}}$
\State $\tilde{\bfmW}^e, \tilde{\bfmW}^d, \tilde{\bfmb}^d \gets \Call{Agglomerate}{\hat{\bfmW}^e, \hat{\bfmW}^d, \hat{\bfmb}^d, \mathcal{M}, \mathcal{M}_n}$
\EndIf
\State Compute loss $\mathcal{L}(\boldsymbol{\mu}_t, \bfmu_{\inmeshlgtwo{n}{t}}; \tilde{\bfmW}^e, \tilde{\bfmW}^d, \tilde{\bfmb}^d, \theta)$
\State Update parameters ${\bfmW}^e, {\bfmW}^d, {\bfmb}^d, \theta$ using gradients $\nabla_{\boldsymbol{W}^e} \mathcal{L},~ \nabla_{\boldsymbol{W}^d} \mathcal{L},~ \nabla_{\boldsymbol{b}^d} \mathcal{L},~ \nabla_{\theta} \mathcal{L}$
\EndFor
\end{algorithmic}
\end{algorithm}

\section{Numerical Simulations}\label{sec:results}
Having introduced the GFN-ROM architecture, providing guarantees on super- and sub-resolution errors, we examine its performance on three challenging benchmarks for MOR. We compare the results of the methodology to intrusive and non-intrusive ROMs: namely, POD-Galerkin (POD-G) and POD with projection (representing the best possible linear method), POD-NN \cite{Hesthaven2018,PichiArtificialNeuralNetwork2023} and GCA-ROM \cite{Pichi2023}. For fairness, we strive to make the MOR baselines for the different models as comparable as possible by adopting similar architectural and hyperparameter choices between methods, as shown in \hyperref[app:hyperparams]{Appendix~\ref*{app:hyperparams}}. Importantly, we do not fine tune the architecture or seek the best possible initialisation for each specific problem. As such, we stress that these results do not represent the best possible performance, but illustrate that GFN-ROM is able to directly handle several challenging benchmarks with no specific fine tuning. As a general rule, we have chosen a latent/reduced dimension $N,L= \lfloor 1.5\times N_\mu \rfloor$, and the small-data regime with a train/test split of 30/70.
 
Each of the problems examined here are selected to showcase a different physical or computational complexity. The first two benchmarks have been previously investigated with GCA-ROM \cite{Pichi2023}, which we use as a baseline for nonlinear reduction -- namely the Graetz problem and an advection-dominated problem. Furthermore, we examine a Stokes flow benchmark with a 7-dimensional parameter space to illustrate the model's capability of dealing with high-dimensional physical and geometric parameter spaces, and the low data regime with sparse samplings of the parameter space. Training data is generated via finite element methods using the RBniCS package \cite{Rozza2024}, which is built on top of the FEniCS package \cite{Alnaes2015, Logg2012}. For each problem, we consider four meshes in order to test the multifidelity performance of our architecture\footnote{The meshes are generated from the largest mesh by using the moving front subsampling algorithm \cite{Lawrence2023}.}. 

The results for single fidelity experiments across our baselines are shown in \autoref{tab:compare_results}.

\begin{table}[h]
    \centering
    \begin{tabular}{c|c|c|c|c|c|c|c|c}
         & POD-G & POD-Proj & POD-NN & GCA-ROM & \multicolumn{4}{c}{GFN-ROM} \\ 
         & Large & Large & Large & Large & Large & Medium & Small & Tiny \\ \hline
        \cellcolor[HTML]{E5E3E3}Graetz & \cellcolor[HTML]{E5E3E3}3.44 & \cellcolor[HTML]{E5E3E3}3.44 & \cellcolor[HTML]{E5E3E3}3.54 & \cellcolor[HTML]{E5E3E3}0.74 & \cellcolor[HTML]{E5E3E3}1.02 & \cellcolor[HTML]{E5E3E3}0.88 & \cellcolor[HTML]{E5E3E3}0.98 & \cellcolor[HTML]{E5E3E3}1.28 \\
        Advection & 27.09 & 25.66 & 30.58 & 4.73 & 4.73 & 4.48 & 7.22 & 12.35 \\
        \cellcolor[HTML]{E5E3E3}Stokes & \cellcolor[HTML]{E5E3E3}2.45 & \cellcolor[HTML]{E5E3E3}2.45 & \cellcolor[HTML]{E5E3E3}2.94 & \cellcolor[HTML]{E5E3E3}4.63 & \cellcolor[HTML]{E5E3E3}4.24 & \cellcolor[HTML]{E5E3E3}4.12 & \cellcolor[HTML]{E5E3E3}4.44 & \cellcolor[HTML]{E5E3E3}5.55 \\ 
    \end{tabular}
    \caption{Mean relative errors (\%) evaluated on large mesh.}
    \label{tab:compare_results}
\end{table}
Since GFN-ROM is capable of being trained on multifidelity data, we investigate its generalisation ability across different resolutions. In \autoref{tab:compare_fid}, we report the change in performance of GFN-ROM when one substitutes half of the training data with cheaper, lower-fidelity data.

\begin{table}[h]
    \centering
    \begin{tabular}{c|c|c|c|c|c|c}
         & \multicolumn{6}{c}{GFN-ROM} \\
         & Large \& & Large \& & Large \& & Medium \& & Medium \& & Small \& \\ 
         & Medium  & Small & Tiny & Small & Tiny & Tiny \\ \hline
        \cellcolor[HTML]{E5E3E3}Graetz & \cellcolor[HTML]{E5E3E3}0.96 \emph{(+0.06)} & \cellcolor[HTML]{E5E3E3}0.98 \emph{(+0.04)} & \cellcolor[HTML]{E5E3E3}1.40 \emph{(-0.37)} & \cellcolor[HTML]{E5E3E3}4.44 \emph{(-3.56)} & \cellcolor[HTML]{E5E3E3}1.03 \emph{(-0.15)} & \cellcolor[HTML]{E5E3E3}1.09 \emph{(-0.11)} \\
        Advection & 5.02 \emph{(-0.30)} & 5.35 \emph{(-0.62)} & 5.63 \emph{(-0.91)} & 5.22 \emph{(-0.73)} & 6.27 \emph{(-1.79)} & 8.77 \emph{(-1.55)} \\
        \cellcolor[HTML]{E5E3E3}Stokes & \cellcolor[HTML]{E5E3E3}3.63 \emph{(+0.61)} & \cellcolor[HTML]{E5E3E3}4.94 \emph{(-0.70)} & \cellcolor[HTML]{E5E3E3}4.44 \emph{(-0.19)} & \cellcolor[HTML]{E5E3E3}4.51 \emph{(-0.39)} & \cellcolor[HTML]{E5E3E3}5.56 \emph{(-1.45)} & \cellcolor[HTML]{E5E3E3}5.65 \emph{(-1.21)} \\
    \end{tabular}
    \caption{Mean relative errors (\%) evaluated on large mesh, and the change when training only with the finer of the two in \autoref{tab:compare_results}. A positive change indicates reduced error and thus better performance.}
    \label{tab:compare_fid}
\end{table}

GFN-ROM is a more lightweight and flexible method compared to existing state-of-the-art ROMs. We show the method's computational efficiency with respect to GCA-ROM in \autoref{tab:compare_efficiency_results}. GFN-ROM offers large savings in terms of training time (up to 63x reduction) and in terms of the number of trainable parameters (up to 25x reduction). All simulations have been performed on a workstation equipped with an Nvidia Quadro RTX 4000 GPU.

\begin{table}[!ht]
    \centering
    \begin{tabular}{cc|c|c|c|c|c}
         &  & \multicolumn{4}{c}{GFN-ROM} \vline & GCA-ROM \\
         &  & Large & Medium & Small & Tiny & Large\\ \hline
        \cellcolor[HTML]{E5E3E3} & \cellcolor[HTML]{E5E3E3}Training time (s) & \cellcolor[HTML]{E5E3E3}119  & \cellcolor[HTML]{E5E3E3}46 & \cellcolor[HTML]{E5E3E3}31 & \cellcolor[HTML]{E5E3E3}26 & \cellcolor[HTML]{E5E3E3}600 \\
        \cellcolor[HTML]{E5E3E3}& \cellcolor[HTML]{E5E3E3}Trainable parameters & \cellcolor[HTML]{E5E3E3}\num{2898761} & \cellcolor[HTML]{E5E3E3}\num{911004} & \cellcolor[HTML]{E5E3E3}\num{311910} & \cellcolor[HTML]{E5E3E3}\num{115821} & \cellcolor[HTML]{E5E3E3}\num{2898795} \\
        \cellcolor[HTML]{E5E3E3}\multirow{-3}{*}{Graetz} & \cellcolor[HTML]{E5E3E3}Mesh nodes & \cellcolor[HTML]{E5E3E3}7205 & \cellcolor[HTML]{E5E3E3}2248 (31\%) & \cellcolor[HTML]{E5E3E3}754 (10\%) & \cellcolor[HTML]{E5E3E3}265 (4\%) & \cellcolor[HTML]{E5E3E3}7205\\ \hline
         \cellcolor[HTML]{FFFFFF}& \cellcolor[HTML]{FFFFFF}Training time (s) & \cellcolor[HTML]{FFFFFF}124 & \cellcolor[HTML]{FFFFFF}53 & \cellcolor[HTML]{FFFFFF}34 & \cellcolor[HTML]{FFFFFF}23 & \cellcolor[HTML]{FFFFFF}408.9 \\
        \cellcolor[HTML]{FFFFFF}& \cellcolor[HTML]{FFFFFF}Trainable parameters & \cellcolor[HTML]{FFFFFF}\num{3538757} & \cellcolor[HTML]{FFFFFF}\num{1110702} & \cellcolor[HTML]{FFFFFF}\num{387298} & \cellcolor[HTML]{FFFFFF}\num{140282} & \cellcolor[HTML]{FFFFFF}\num{3538791}\\
        \cellcolor[HTML]{FFFFFF}\multirow{-3}{*}{Advection}& \cellcolor[HTML]{FFFFFF}Mesh nodes & \cellcolor[HTML]{FFFFFF}8801 & \cellcolor[HTML]{FFFFFF}2746 (31\%) & \cellcolor[HTML]{FFFFFF}942 (11\%) & \cellcolor[HTML]{FFFFFF}326 (4\%) & \cellcolor[HTML]{FFFFFF}8801 \\ \hline
        \cellcolor[HTML]{E5E3E3} & \cellcolor[HTML]{E5E3E3}Training time (s) & \cellcolor[HTML]{E5E3E3}248 & \cellcolor[HTML]{E5E3E3}90 & \cellcolor[HTML]{E5E3E3}45 & \cellcolor[HTML]{E5E3E3}31 & \cellcolor[HTML]{E5E3E3}1949 \\
        \cellcolor[HTML]{E5E3E3}& \cellcolor[HTML]{E5E3E3}Trainable parameters & \cellcolor[HTML]{E5E3E3}\num{2827589} & \cellcolor[HTML]{E5E3E3}\num{905596} & \cellcolor[HTML]{E5E3E3}\num{316126} &  \cellcolor[HTML]{E5E3E3}\num{118032} & \cellcolor[HTML]{E5E3E3}\num{2827623} \\
        \cellcolor[HTML]{E5E3E3}\multirow{-3}{*}{Stokes}& \cellcolor[HTML]{E5E3E3}Mesh nodes & \cellcolor[HTML]{E5E3E3}7019 & \cellcolor[HTML]{E5E3E3}2226 (32\%) & \cellcolor[HTML]{E5E3E3}756 (11\%) & \cellcolor[HTML]{E5E3E3}262 (4\%) & \cellcolor[HTML]{E5E3E3}7019\\ \hline
    \end{tabular}
    \caption{Computational efficiency of GFN-ROM in comparison to GCA-ROM, and its scaling w.r.t.\ the mesh size.}
\label{tab:compare_efficiency_results}
\end{table}

\subsection{Graetz problem}
\begin{figure}[h]
\centering
\begin{subfigure}[c]{.42\textwidth}
\centering
\includegraphics{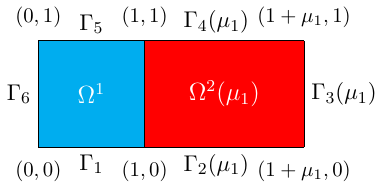}
\caption{}
\label{fig:graetz_domain}
\end{subfigure}
\hfill
\begin{subfigure}[c]{.56\textwidth}
    \centering
\includegraphics{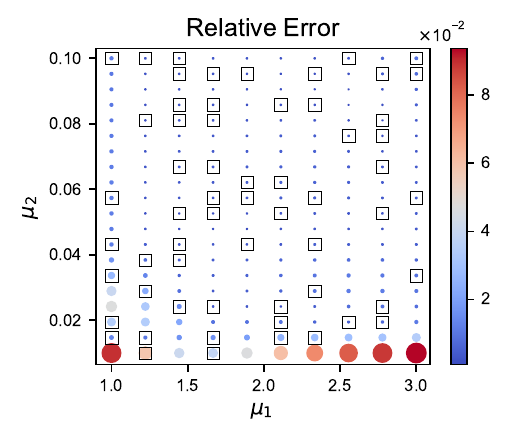}
\caption{}
\label{fig:graetz_error_field}
\end{subfigure}
\caption{(A) The parameterised domain for the Graetz problem. (B) Mean relative errors for GFN-ROM trained on the large mesh for the Graetz problem. Squares indicate the parameter realisations seen during training.}
\end{figure}

The Graetz problem represents a common benchmark in MOR, combining a forced heat diffusion with horizontal heat advection in a parametrised geometry. Specifically, we study a steady-state Graetz problem on the parametrised geometry $\Omega(\mu_1) = \Omega^1 \cup \Omega^2(\mu_1) \in \mathbb{R}^2$, with $\Omega^1 = [0,1]\times[0,1]$ and $\Omega^2(\mu_1)=[1,1+\mu_1]\times[0,1]$  illustrated in \autoref{fig:graetz_domain}, and discretised by four meshes as summarised in Table \ref{tab:compare_efficiency_results}.

Homogeneous Dirichlet conditions are imposed on $\Gamma_D = \Gamma_1 \cup \Gamma_5 \cup \Gamma_6$, whilst a constant Dirichlet boundary condition is imposed on $\Gamma_G =\Gamma_2(\mu_1) \cup \Gamma_4(\mu_1)$ as $u(\boldsymbol{\mu})|_{\Gamma_G}=1$. The remaining boundary $\Gamma_3(\mu_1)$ is given a homogeneous von Neumann boundary condition. The weak formulation of the Graetz problem reads as follows: given $\boldsymbol{\mu} \in \mathbb{P}$, find $u(\boldsymbol{\mu}) \in \mathbb{U}(\mu_1)$
s.t.\
\bea
\mu_2 \int_{\Omega(\mu_1)} \nabla u \cdot \nabla v\, d\Omega + \int_{\Omega(\mu_1)} y(1-y) \frac{\partial u}{\partial x} v\, d\Omega = 0, \quad \forall v \in \mathbb{V},
\eea
where $\mathbb{U}(\mu_1) = \{v \in H^1(\Omega(\mu_1)): v |_{\Gamma_D}=0, ~ v |_{\Gamma_G}=1 \}$ and $\mathbb{V}(\mu_1) = \{v \in H^1(\Omega(\mu_1)): v |_{\Gamma_D \cup \Gamma_G} =0\}$. Note that due to the parametrised domain, these function spaces are dependent on the geometric parameter.
This problem is parametrised by $\boldsymbol{\mu} = (\mu_1, \mu_2)$, where $\mu_1$ is a geometrical parameter governing the length of $\Omega^2$, and $\mu_2$ is a physical parameter, namely the diffusivity coefficient. We consider a dataset with 200 snapshots computed on a uniform grid of $\mathbb{P}$, with 10 and 20 equispaced values for $\mu_1 \in [1,3]$ and $\mu_2 \in [0.01, 0.1]$, respectively.

Training and testing the performance on GFN-ROM on the largest mesh, we achieve a mean relative error of 1.02\% over the full dataset. As illustrated in \autoref{fig:graetz_error_field}, we can see GFN-ROM is well able to capture the overall solution. Note that this performance is very good, particularly considering only 60 parameter realisations are seen for training, and the model learns a latent representation of size $L=3$. Interestingly, for the Graetz problem we see that errors are much higher for small values of $\mu_2$ i.e.\ in the low diffusivity region, where the advection dominates. Moreover, few training samples are representing this regime, making it more difficult for GFN-ROM to learn how to model the region.

Nonetheless, the overall performance on GFN-ROM is extremely good on this benchmark, achieving over a 3x reduction in mean relative error compared to linear methods (see \autoref{tab:compare_results}), and comparable performance to GCA-ROM. 

\begin{figure}[h]
    \centering
    \begin{subfigure}[b]{.49\textwidth}
    \centering
\includegraphics{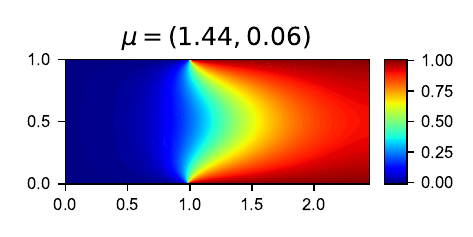}
     \caption{Training on large mesh.}
    \end{subfigure}\hfill
    \begin{subfigure}[b]{.49\textwidth}
    \centering
\includegraphics{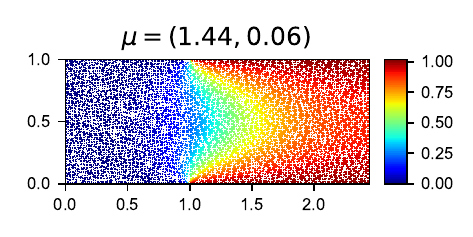}
     \caption{Training on medium mesh.}
    \end{subfigure}
\begin{subfigure}[b]{.49\textwidth}
    \centering
\includegraphics{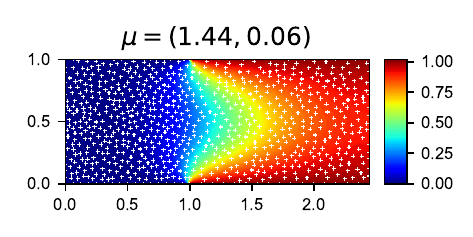}
     \caption{Training on small mesh.}
    \end{subfigure}\hfill
    \begin{subfigure}[b]{.49\textwidth}
    \centering
\includegraphics{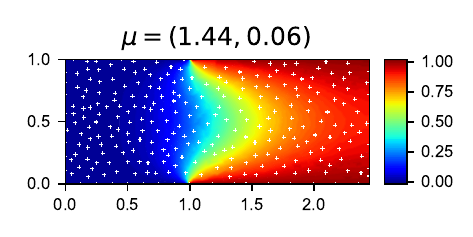}
     \caption{Training on tiny mesh.}
    \end{subfigure}
    \caption{GFN-ROM solution fields for the Graetz problem when trained on (A) large mesh, (B) medium mesh, (C) small mesh, and (D) tiny mesh. The mesh points are shown by white pluses.}
    \label{fig:graetz_examples}
\end{figure}

Moreover, we can further improve the performance of GFN-ROM by leveraging its resolution invariance. In fact, when training on a coarser representation, e.g.\ the medium mesh, the mean relative error is reduced to 0.88\%, illustrating the new opportunities arising from GFN-ROM's ability to perform super-resolution. Even when training on the coarsest (tiny) discretisation, the mean relative error increases only slightly by 0.26\% to 1.28\% despite this representing only 4\% of the training data. In \autoref{fig:graetz_examples}, we show the predictions of the four different GFM-ROM models trained on the four different meshes for a parameter realisation $\boldsymbol{\mu}=(1.44,0.06)$ not seen during training. The model predictions are qualitatively indistinguishable, showing that GFN-ROM is capable of extracting enough information even on cheaper computational data without incurring any penalty on performance. Even further generalisations are possible thanks to the  possibility of training on mixed-fidelity data. GFN-ROM shows excellent generalisation properties and performance is minorly impacted when substituting high-fidelity data for lower-fidelity data, as seen in \autoref{tab:compare_fid}. For example, the performance even improves from 1.02\% when training solely on the large mesh to 0.96\% when substituting half of the large mesh dataset with cheaper data from the medium mesh (representing a 35\% reduction in training data).

\subsection{Advection-dominated problem}
\begin{figure}[th]
\centering
\begin{subfigure}[c]{.42\textwidth}
\centering
\includegraphics{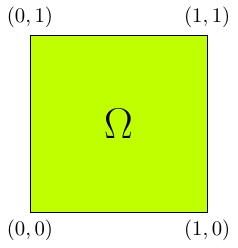}
\caption{}
\label{fig:adv_domain}
\end{subfigure}
\hfill
\begin{subfigure}[c]{.56\textwidth}
    \centering
\includegraphics{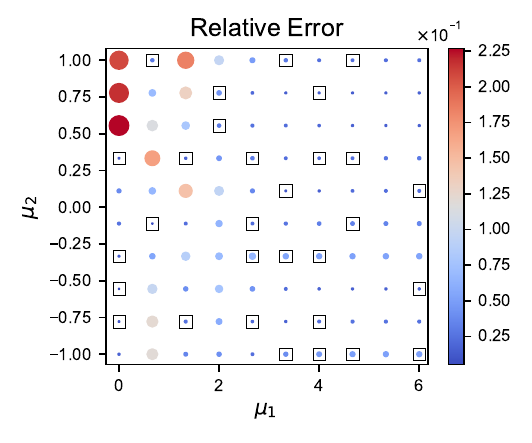}
\caption{}
\label{fig:advection_error_field}
\end{subfigure}
\caption{(A) The parameterised domain for the advection problem. (B) Mean relative errors for GFN-ROM trained on the large mesh for the advection problem. Squares indicate the parameter realisations seen during training.}
\end{figure}

Advection-dominated problems, due to the slow decay of the Kolmogorov $n$-width \cite{Greif2019, Ohlberger2015}, represent an important benchmarking case in the ROM context to compare linear and nonlinear reduction strategies. We consider a fixed domain $\Omega=[0,1]\times[0,1]\in\mathbb{R}^2$, illustrated in \autoref{fig:adv_domain} with homogeneous Dirichlet boundary conditions. We generate four meshes summarised in Table \ref{tab:compare_efficiency_results}.

The weak formulation of the advection problem reads as follows: 
given $\boldsymbol{\mu} \in \mathbb{P}$, find $u(\boldsymbol{\mu}) \in \mathbb{V}$ s.t.\
\bea
{\operatorname{D}(\mu_1)} \int_\Omega \nabla u \cdot \nabla v\, d\Omega + \int_\Omega (\boldsymbol{\beta} \cdot \nabla u) v\, d\Omega = \int_\Omega v\, d\Omega, \quad \forall v \in \mathbb{V}= H^1_0(\Omega).
\eea
A total of 100 snapshots are generated for this problem on a uniform grid of $\mathbb{P}$, with 10 and 10 equispaced values for $\mu_1 \in [0,6]$ and $\mu_2 \in [-1, 1]$, respectively. The first parameter $\mu_1$ determines the diffusion coefficient via $\operatorname{D}(\mu_1)=10^{-\mu_1}$ and the second parameter $\mu_2$ determines the strength and sign of the transport $\boldsymbol{\beta}(\mu_2)=(\mu_2,\mu_2)$.

As expected, POD-based techniques perform poorly due to their linear nature on this benchmark, with POD-G only achieving a mean relative error of 27.09\% over the full dataset. GFN-ROM's nonlinear nature allows it to achieve a far superior mean relative error of only 4.73\% when training on the large mesh, which matches GCA-ROM's performance. We show the relative errors over the parameter space in \autoref{fig:advection_error_field}. We stress that this performance is possible even in the low data regime with only 30 training examples and a small latent representation of size $L=3$. Contrarily to the Graetz problem, it can be seen that poorer performances occur for lower Péclet numbers ($\operatorname{Pe}\approx \lvert \frac{\mu_2}{10^{-\mu_1}} \rvert$) compared to higher Péclet numbers i.e.\ the model performs better when transport dominates. This is likely due to the fact that the model training data is skewed towards examples where the advection in the opposite direction dominates, making it less adept at predicting solutions where diffusion has a stronger impact and the $\mu_2$ is positive.

\begin{figure}[h]
    \centering
    \begin{subfigure}[b]{.49\textwidth}
    \includegraphics{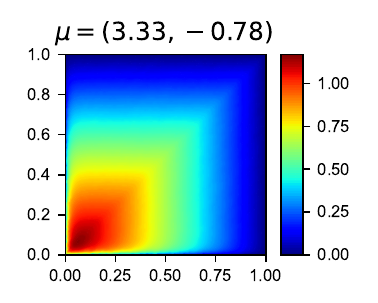}
     \caption{Training on large mesh.}
    \end{subfigure}
    \hfill
    \begin{subfigure}[b]{.49\textwidth}
    \includegraphics{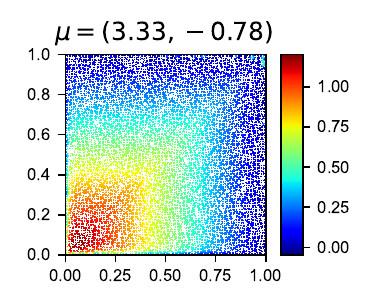}
     \caption{Training on medium mesh.}
    \end{subfigure}
    \begin{subfigure}[b]{.49\textwidth}
    \includegraphics{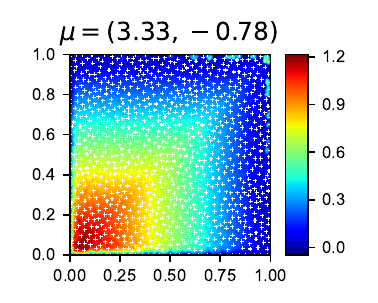}
     \caption{Training on small mesh.}
    \end{subfigure}
    \hfill
    \begin{subfigure}[b]{.49
\textwidth}
    \includegraphics{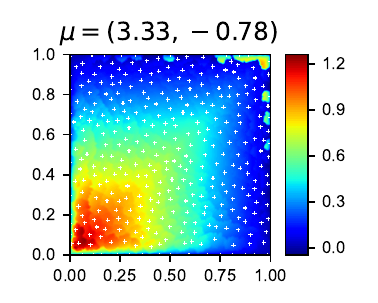}
     \caption{Training on tiny mesh.}
    \end{subfigure}
    \caption{GFN-ROM solution fields for the advection problem when trained on (A) large mesh, (B) medium mesh, (C) small mesh, and (D) tiny mesh. The mesh points are shown by white pluses.}
    \label{fig:advection_examples}
\end{figure}

On this benchmark, GFN-ROM's resolution invariant nature allows for training on cheaper meshes, and outperforming  GCA-ROM, achieving a mean relative error of 4.48\% when training on the medium mesh (69\% reduction in training data). In \autoref{fig:advection_examples}, we show the predictions of the GFN-ROM models trained on the four meshes for a parameter $\boldsymbol{\mu}=(3.33,-0.78)$ not seen during training. Qualitatively, we can see that all models but tiny are capable of learning well-enough the solution dynamics. The tradeoff to consider is of course given by the performance of large models and the computational efficiency of the smaller ones. Similarly to the Graetz benchmark, we also observed good generalisation for the multifidelity cases, as shown in \autoref{tab:compare_fid}.

\subsection{Stokes problem}
\begin{figure}[!t]
\begin{minipage}[b]{\textwidth}
\centering
\begin{tikzpicture}[scale=0.7]
            \node[left] at (-1, 3) {\large{$(0, \mu_0+ \mu_2 + \mu_3)$}}; 
            \node[right] at (1, 3) {\large{$(\mu_1, \mu_0+ \mu_2 + \mu_3)$}}; 
            \node[right] at (3, 1) {\large{$(\mu_1 + \mu_4, \mu_0+ \mu_2 + \mu_4\tan(\mu_5))$}}; 
            \node[right] at (3, -1) {\large{$(\mu_1 + \mu_4, \mu_2 + \mu_4\tan(\mu_5))$}}; 
            \node[left] at (-1, 1) {\large{$(0, \mu_0+\mu_2)$}}; 
            \node[left] at (-1, -1) {\large{$(0, \mu_2)$}}; 
            \node[left] at (-1, -3) {\large{$(0, 0)$}}; 
            \node[right] at (1, -3) {\large{$(\mu_1, 0)$}}; 
            \fill[green!55!blue] (0,-1) rectangle (3,1);
            \fill[green!55!blue] (-1,-3) rectangle (1,3);
            \node[white] at (1,0) {\LARGE{$\Omega(\boldsymbol{\mu}_{\text{geo}})$}};
            \draw (-1,-3) -- (-1,3) node[pos=0.5,left] {\Large{$\Gamma_{\text{w}}$}};
            \draw (-1,3) -- (1,3) node[white,pos=0.5,below] {\Large{$\Gamma_{\text{N}}$}};
            \draw (-1,-3) -- (1,-3) node[white,pos=0.5,above] {\Large{$\Gamma_{\text{N}}$}};
            \draw (1,-1) -- (3,-1) node[pos=0.5,below] {\Large{$\Gamma_{\text{w}}$}};
            \draw (1,1) -- (3,1) node[pos=0.5,above] {\Large{$\Gamma_{\text{w}}$}};            
            \draw (1,-3) -- (1,-1);
            \draw (1,1) -- (1,3);
            \draw (3,-1) -- (3,1) node[pos=0.5,right] {\Large{$\Gamma_{\text{in}}$}};
            \end{tikzpicture}
\captionof{figure}{The domain for the Stokes problem.}
\label{fig:stokes_domain}
\end{minipage}
\end{figure}
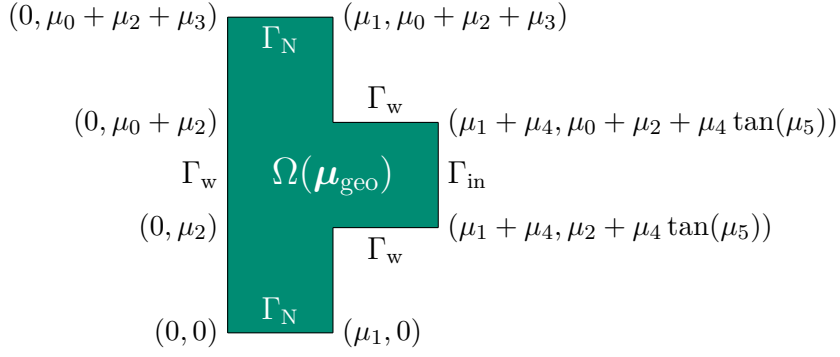

We now consider a more complex benchmark in computational fluid dynamics, modelling a steady-state Stokes flow problem with high-dimensional physical/geometric parameter space and sparse sampling of it. This setting makes the problem particularly challenging to learn, and an important test case for the ability of a ROM to deal with the low-data regime. For simplicity, we consider the case of a viscosity of fixed value $1$, but we allow for a variable forcing term $\bfmf = (0, \mu_6)$, governed by a physical parameter $\mu_6$. We also consider a total of 6 geometrical parameters $\boldsymbol{\mu}_{\text{geo}} = (\mu_0, \mu_1, \mu_2, \mu_3, \mu_4, \mu_5)$, to deform the domain depicted in \autoref{fig:stokes_domain}. We generate four meshes summarised in \autoref{tab:compare_efficiency_results}.
As boundary conditions, the inlet $\Gamma_{\text{in}}$ is given a Dirichlet boundary condition of $\bfmu_{in} = [4(y-1)(2-y), 0]$, the walls $\Gamma_{\text{w}}$ are given zero Dirichlet boundary conditions and the outlet $\Gamma_{\text{N}}$ is left as a homogeneous Neumann boundary condition. A weak formulation for this problem reads: for a given parameter $\boldsymbol{\mu} \in\mathbb{P}$, find $\boldsymbol{u}(\boldsymbol{\mu}) \in\mathbb{U}$, $p \in\mathbb{M}$ s.t.\ 

\bea
\begin{cases}
    \int_{\Omega} \nabla \boldsymbol{u} : \nabla \boldsymbol{v} \, d\Omega - \int_{\Omega} p \nabla \cdot \boldsymbol{v} \, d\Omega = \int_{\Omega} \boldsymbol{f} \cdot \boldsymbol{v} \, d\Omega, \quad \forall \boldsymbol{v} \in\mathbb{V},  \\
    \int_{\Omega} q \nabla \cdot \boldsymbol{u} \, d\Omega = 0, \quad \forall q \in\mathbb{M},
\end{cases}
\eea
where the parameter dependent function spaces for the velocity are defined as $\mathbb{V}(\boldsymbol{\mu}_{\text{geo}}) = \{\bfmv \in (H^1(\Omega(\boldsymbol{\mu}_{\text{geo}}))^2: \bfmv|_{\Gamma_{w}\cup\Gamma_{in}}=\boldsymbol{0} \}$ and $\mathbb{U}(\boldsymbol{\mu}_{\text{geo}}) = \{\bfmv \in (H^1(\Omega(\boldsymbol{\mu}_{\text{geo}})))^2: \bfmv|_{\Gamma_{w}}=\boldsymbol{0}, \bfmv|_{\Gamma_{in}} = \bfmu_{in}\}$, and the function space for the pressure is defined as $\mathbb{M}(\boldsymbol{\mu}_{\text{geo}}) = L^2(\Omega(\boldsymbol{\mu}_{\text{geo}}))$. The problem uses a mixed finite element discretisation with the velocity and pressure as solution variables, meaning that the inf-sup condition is necessary for the well posedness of this problem and as a result the supremiser operator $T^{\mu}: \mathbb{M} \rightarrow \mathbb{V}$ is used (see \cite{BallarinSupremizerStabilizationPOD2015} for more details).

We consider a sparse sampling with solutions computed only for the following parameter choices: $\mu_0, \mu_1, \mu_2, \mu_3, \mu_4 \in \{0.5,1.5 \}, \mu_5 \in \{-\frac{\pi}{6}, \frac{\pi}{6} \}, \mu_6 \in \{-10,-8,-6,-4,-2,0,2,4,6,8,10\}$, leading to a total of 704 high-fidelity solutions, and the field of interest is the magnitude of the the velocity $|\bfmu|$.

\begin{figure}[h]
    \centering
    \includegraphics{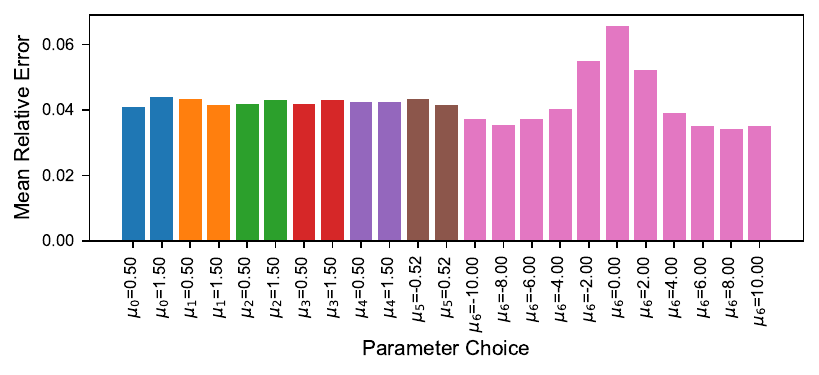}
    \caption{Histogram of the mean relative error for GFN-ROM trained on the large mesh for the Stokes problem.}
    \label{fig:stokes_errors}
\end{figure}

Training on the large mesh, GFN-ROM achieves a mean relative error of 4.24\%, outperforming GCA-ROM. We stress that this performance is possible whilst remaining in the low data regime, undersampling the parameter space. Despite this, the model is well able to capture the solution. As shown in \autoref{fig:stokes_errors}, model performance is found to be very similar on average for the possible geometrical parameters, while, as expected, the results are more sensitive to the physical parameter $\mu_6$, with slightly worse approximation for the case $\mu_6=0$ where the forcing term vanishes.

\begin{figure}[htb]
    \centering
    \begin{minipage}{.3\textwidth}
\includegraphics{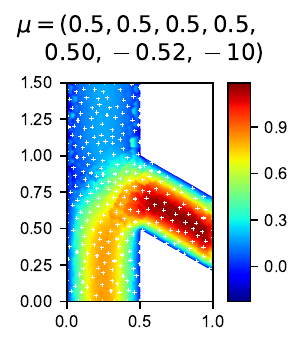}
    \end{minipage}
    \hfill
        \begin{minipage}{.3\textwidth}
\includegraphics{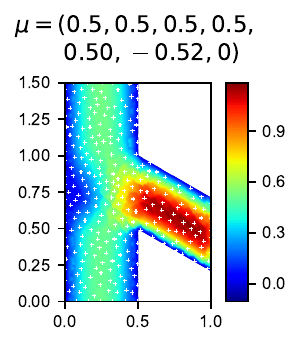}
    \end{minipage}
        \hfill
        \begin{minipage}{.3\textwidth}
\includegraphics{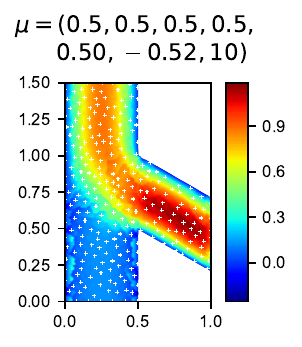}
    \end{minipage}
    \caption{Super-resolution via GFN-ROM when training on the tiny mesh and evaluating on large mesh for variations in the physical parameter $\mu_6$ showcasing downward, none and upward forcings.}
    \label{fig:stokes_fields_1}
\end{figure}

\begin{figure}[htb]
    \centering
    \begin{minipage}{.3\textwidth}
\includegraphics{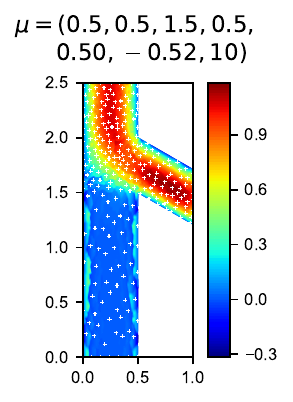}
    \end{minipage}
    \hfill
        \begin{minipage}{.3\textwidth}
\includegraphics{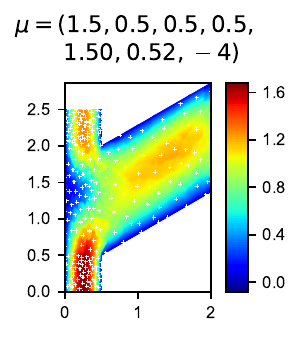}
    \end{minipage}
        \hfill
        \begin{minipage}{.3\textwidth}
\includegraphics{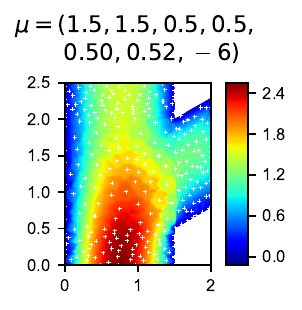}
    \end{minipage}
    \caption{Super-resolution via GFN-ROM when training on the tiny mesh and evaluating on large mesh for variations in the geometrical parameters.}
    \label{fig:stokes_fields_2}
\end{figure}

GFN-ROM shows even further performance improvements by training on cheaper computational meshes, achieving a mean relative error of 4.12\% for the medium mesh (representing a 68\% reduction in the number of nodes). Even when training on the coarsest discretisation with just 4\% of the data, performance only increases to 5.55\% and is still able to capture the behaviour with respect to the physical and geometric parameters not seen during training, respectively in \autoref{fig:stokes_fields_1} and \autoref{fig:stokes_fields_2}.

Once again, the complex parametric variability and the scarce amount of information (training on 211 parameter realisations for the 7-dimensional parameter space) do not compromise the generalisation capability of the methodology. 

We also see excellent results for multifidelity training. In general, performances are only minorly affected when substituting half of the data for cheaper data, underlining the possibility of training on cheaper data whilst retaining excellent performance. It is even possible to achieve a better performance of 3.63\% by replacing 50\% of data on the large mesh with cheaper data on the medium mesh (representing a 34\% reduction in training data), as detailed in \autoref{tab:compare_fid}.

\section{Conclusion}\label{sec:conc}

In this work, we presented a novel  graph-based resolution-invariant ROM, as a powerful tool to investigate multifidelity applications in the MOR context. We showed that it is possible to apply feedforward networks for graphical data by means of the graph feedforward network, allowing for the extension of many single-fidelity architectures to the multifidelity setting. Furthermore, we showed that this extension comes with provable guarantees on the performance in terms of error bounds.

Beyond such theoretical guarantees, we have also demonstrated that our architecture performs well in practice on three challenging benchmarks involving parametrised domains, advection-dominated problems and high-dimensional parameter spaces. Across all benchmarks, GFN-ROM achieved satisfactory performances, especially considering we worked within a low-data regime. Indeed, like GCA-ROM, our architecture has been trained on only 30\% of the dataset, and has been exploited to learn a small latent representation.
GFN-ROM achieved comparable or better performance than GCA-ROM whilst being more flexible, interpretable and lightweight, and applicable for multifidelity data. 

We showed that GFN-ROM can perform super-resolution, and demonstrated that within the proposed architecture, it is possible to substitute high-fidelity data for lower-fidelity data without incurring performance penalties. In practice, we showed that often performance can even be improved when training on cheaper computational data. Our work thus highlighted the importance of multifidelity ROMs, showing that computational
savings can be made by training on smaller meshes without deterioration in performance.

Future work will investigate further possibilities with graph feedforward networks, including tackling time-dependent problems, problems with spatially dependent parameters, enforcing more smoothness in GFN-ROM predictions, forcing greater dependence on local connections and automating finding the best training mesh (adapting the number of model parameters during training). We also plan to investigate multimodal modelling by using information from the Fourier domain, and the generalisation of convolutional neural operators to unstructured grids.

\section*{Acknowledgment}
FP acknowledges the “GO for IT” program within the CRUI fund for the project “Reduced order method for nonlinear PDEs enhanced by machine learning”. This work has been conducted within the research activities of the consortium iNEST (Interconnected North-East Innovation Ecosystem), Piano Nazionale di Ripresa e Resilienza (PNRR) – Missione 4 Componente 2, Investimento 1.5 – D.D. 1058 23/06/2022, ECS00000043, supported by the European Union's NextGenerationEU program.

\bibliographystyle{abbrv}
\bibliography{main}

\newpage
\appendix

\section{Proofs}

\subsection{Extension to multiple layers}
\label{app:mult_layers}
For ease of notation we define each of the outputs of the $\operatorname{GFN}^{\mathcal{M}_o \to \mathcal{M}_n}$ function as 
\nbea
\operatorname{GFN}_{W^e}^{\mathcal{M}_o \to \mathcal{M}_n}(\bfmW^e), ~\operatorname{GFN}_{W^d}^{\mathcal{M}_o \to \mathcal{M}_n}(\bfmW^d), ~\operatorname{GFN}_{b^d}^{\mathcal{M}_o \to \mathcal{M}_n}(\bfmb^d) = \operatorname{GFN}^{\mathcal{M}_o \to \mathcal{M}_n}(\bfmW^e, \bfmW^d, \bfmb^d).
\neea
The single-layer GFN autoencoder given in \autoref{eq:enc_Mn} and \autoref{eq:dec_Mn} can then be rewritten as
\nbea
\operatorname{enc}^{\mathcal{M}_o \to \mathcal{M}_n}(\bfmu_{\inmeshlg{n}})_i &=& \sigma \left( \sum_{j_{\inmesh{n}}=1}^{N_n} \operatorname{GFN}_{W^e}^{\mathcal{M}_o \to \mathcal{M}_n}(\bfmW^e)_{ij_{\inmesh{n}}} u_{j_{\inmesh{n}}} + b^e_i \right), \qquad \forall i=1,\cdots,L,\\
\operatorname{dec}^{\mathcal{M}_o \to \mathcal{M}_n}(\bfmz)_{i_{\inmesh{n}}} &=& \sum_{j=1}^{L} \operatorname{GFN}_{W^d}^{\mathcal{M}_o \to \mathcal{M}_n}(\bfmW^d)_{i_{\inmesh{n}}j} z_j + b^d_{i_{\inmesh{n}}}, \qquad \forall i_{\inmeshlg{n}}=1,\cdots,N_n.
\neea
This can be extended to $Q+1$ hidden layers as
\nbea
\operatorname{enc}_0^{\mathcal{M}_o \to \mathcal{M}_n}(\bfmx)_i &=& \sigma \left( \sum_{j_{\inmesh{n}}=1}^{N_n} \operatorname{GFN}_{W^e}^{\mathcal{M}_o \to \mathcal{M}_n}(\bfmW^e)_{ij_{\inmesh{n}}} x_{j_{\inmesh{n}}} + b^e_i \right) , \nonumber \\
\label{eq:gfn_ae1}
&& \forall {i=1,\cdots,L_1}, \\
\operatorname{enc}_{P\cdots 1}(\bfmx)_i &=& \sigma \left( \sum_{l_P=1}^{L_P} W^{(P)}_{i,l_P} \sigma \left( \cdots \sigma \left( \sum_{l_1=1}^{L_1} W^{(1)}_{l_2, l_1} x_{l_1} + b^{(1)}_{l_2} \right) + \cdots \right) + b^{(P)}_i \right), \nonumber \\
&& \forall {i=1,\cdots,L_{P+1}}, \\
\operatorname{dec}_{Q\cdots P+1}(\bfmx)_i &=& \sigma \left( \sum_{l_Q=1}^{L_Q} W^{(Q)}_{i, l_Q} \sigma \left( \cdots \sigma \left( \sum_{l_{P+1}=1}^{L_{P+1}} W^{(P+1)}_{l_{P+2}, l_{P+1}} x_{l_{P+1}} + b^{(P+1)}_{l_{P+2}} \right) + \cdots \right) + b^{(Q)}_i \right), \nonumber \\
&& \forall {i=1,\cdots,L_{Q+1}}, \\
\operatorname{dec}_{Q+1}^{\mathcal{M}_o \to \mathcal{M}_n}(\bfmx)_{i_{\inmesh{n}}} &=& \sum_{j=1}^{L_{Q+1}} \operatorname{GFN}_{W^d}^{\mathcal{M}_o \to \mathcal{M}_n}(\bfmW^d)_{i_{\inmesh{n}}j} x_j + \operatorname{GFN}_{b^d}^{\mathcal{M}_o \to \mathcal{M}_n}(\bfmb^d)_{i_{\inmesh{n}}}, \nonumber \\
\label{eq:gfn_ae4}
&& \forall i_{\inmeshlg{n}}=1,\cdots,N_n,
\neea
where $0 \leq P\leq Q$, $\mathcal{M}_n$ is a new mesh, $L_{i}$ represents the output sizes for the $(i+1)$-th network layer with $L_{P+1}=L$ and $L_{Q+2}=N_n$, $\operatorname{enc}^{\mathcal{M}_o \to \mathcal{M}_n}(\bfmx)_i = \operatorname{enc}_{P\cdots 1}(\operatorname{enc}^{\mathcal{M}_o \to \mathcal{M}_n}_{0}(\bfmx))_i, \forall i=1,\cdots,L$ is the encoded latent representation of an input $\bfmx \in \mathbb{R}^{N_n}$ on mesh $\mathcal{M}_n$ and $\operatorname{dec}^{\mathcal{M}_o \to \mathcal{M}_n}(\bfmz)_{i_{\inmesh{n}}}=\operatorname{dec}^{\mathcal{M}_o \to \mathcal{M}_n}_{Q+1}(\operatorname{dec}_{Q\cdots P+1}(\bfmz))_{i_{\inmesh{n}}}, \forall i_{\inmeshlg{n}}=1,\cdots,N_n$ is the reconstructed solution on $\mathcal{M}_n$ from a given latent representation $\bfmz \in \mathbb{R}^{L}$.

\subsection{Simplified expansive equations}
\label{subsec:exp_proof}
We consider the GFN transforms $\mathcal{M}_o \to \mathcal{M}_n$ and show that they can be simplified under the expansive condition given in \autoref{eq:exp}.

First of all, we show that the condition implies
\nbea
\forall i_{\inmeshlg{n}}, \quad \{ k_{\inmeshlg{o}} \text{ s.t. } k_{\inmeshlg{o}} \eitherarrow i_{\inmeshlg{n}} \} = \{ k_{\inmeshlg{o}} \text{ s.t. } k_{\inmeshlg{o}} \leftarrow i_{\inmeshlg{n}} \},
\neea
which is a set containing a unique element $k_{\inmeshlg{o}}$ defined as $k_{\inmeshlg{o}} \leftarrow i_{\inmeshlg{n}}$.

\emph{Proof:}
We have to show that $\{ k_{\inmeshlg{o}} \text{ s.t. } k_{\inmeshlg{o}} \eitherarrow i_{\inmeshlg{n}} \} = \{ k_{\inmeshlg{o}} \text{ s.t. } k_{\inmeshlg{o}} \leftarrow i_{\inmeshlg{n}} \} \cup \{ k_{\inmeshlg{o}} \text{ s.t. } k_{\inmeshlg{o}} \rightarrow i_{\inmeshlg{n}} \} = \{ k_{\inmeshlg{o}} \text{ s.t. } k_{\inmeshlg{o}} \leftarrow i_{\inmeshlg{n}} \}$. If $\{ k_{\inmeshlg{o}} \text{ s.t. } k_{\inmeshlg{o}} \rightarrow i_{\inmeshlg{n}} \} = \emptyset$, this is trivially the case. Therefore, in the following we consider only the case where it is not.

Trivially, $|\{ k_{\inmeshlg{o}} \text{ s.t. } k_{\inmeshlg{o}} \leftarrow i_{\inmeshlg{n}} \}|=1$ since a node $i_{\inmeshlg{n}}$ can only have one nearest neighbour. In fact, additionally we can show it is also the case that $|\{ k_{\inmeshlg{o}} \text{ s.t. } k_{\inmeshlg{o}} \rightarrow i_{\inmeshlg{n}} \}|=1$, which can be proven by contradiction. Suppose that $|\{ k_{\inmeshlg{o}} \text{ s.t. } k_{\inmeshlg{o}} \rightarrow i_{\inmeshlg{n}} \}|=R \neq 1$. This means that there exist multiple distinct nodes $\left\{k^r_{\inmeshlg{o}}\right\}_{r=1}^R$, each satisfying $k^r_{\inmeshlg{o}} \rightarrow i_{\inmeshlg{n}}$. Via the expansive condition in \autoref{eq:exp}, it follows that each of these nodes must also satisfy $k^r_{\inmeshlg{o}} \leftarrow i_{\inmeshlg{n}}$. This is not possible since $i_{\inmeshlg{n}}$ can only have a single nearest neighbour, meaning therefore $|\{ k_{\inmeshlg{o}} \text{ s.t. } k_{\inmeshlg{o}} \rightarrow i_{\inmeshlg{n}} \}|=1$.

We note that the single node in the set $|\{ k_{\inmeshlg{o}} \text{ s.t. } k_{\inmeshlg{o}} \rightarrow i_{\inmeshlg{n}} \}|$ must also satisfy $k_{\inmeshlg{o}} \leftarrow i_{\inmeshlg{n}}$ due to \autoref{eq:exp}. Therefore, it is in fact the exact same node as the node in the set $|\{ k_{\inmeshlg{o}} \text{ s.t. } k_{\inmeshlg{o}} \leftarrow i_{\inmeshlg{n}} \}|$, thus concluding the proof.

Secondly, we show that for a given $k_{\inmeshlg{o}}$
\nbea
\lvert \{ h_{\inmeshlg{n}} \text{ s.t. } k_{\inmeshlg{o}} \eitherarrow h_{\inmeshlg{n}} \}\rvert = \lvert \{ h_{\inmeshlg{n}} \text{ s.t. } k_{\inmeshlg{o}} \leftarrow h_{\inmeshlg{n}} \}\rvert.
\neea

\emph{Proof:}
Due to \autoref{eq:exp}, it is the case that $\{ h_{\inmeshlg{n}} \text{ s.t. } k_{\inmeshlg{o}} \rightarrow h_{\inmeshlg{n}} \} \subseteq \{ h_{\inmeshlg{n}} \text{ s.t. } k_{\inmeshlg{o}} \leftarrow h_{\inmeshlg{n}} \}$, concluding the proof.

The updates can therefore be simplified to
\nbea
\tilde{W}^e_{ij_{\inmesh{n}}} &=& \frac{1}{|\{h_{\inmeshlg{n}} \text{ s.t. } k_{\inmeshlg{o}} \leftarrow h_{\inmeshlg{n}} \}|} {W}^e_{ik_{\inmesh{o}}}, \qquad \text{ where}\quad k_{\inmeshlg{o}} \leftarrow j_{\inmeshlg{n}}, \\
\tilde{W}^d_{i_{\inmesh{n}}j} &=& {W}^d_{k_{\inmesh{o}}j}, \qquad \text{ where}\quad k_{\inmeshlg{o}} \leftarrow i_{\inmeshlg{n}}, \\
\tilde{b}^d_{i_{\inmesh{n}}} &=& {b}^d_{k_{\inmesh{o}}},  \qquad \text{ where}\quad k_{\inmeshlg{o}} \leftarrow i_{\inmeshlg{n}},
\neea
which therefore concludes the proof.

\subsection{Simplified agglomerative equations}
\label{subsec:agg_proof}
We consider the GFN transforms $\mathcal{M}_o \to \mathcal{M}_n$ and show that they can be simplified under the agglomerative condition given in \autoref{eq:agg}.

First of all, we show that the condition implies
\nbea
\forall j_{\inmeshlg{n}}, \quad \{ k_{\inmeshlg{o}} \text{ s.t. } k_{\inmeshlg{o}} \eitherarrow j_{\inmeshlg{n}} \} = \{ k_{\inmeshlg{o}} \text{ s.t. } k_{\inmeshlg{o}} \rightarrow j_{\inmeshlg{n}} \}.
\neea

\emph{Proof:}
Due to \autoref{eq:agg}, it is the case that $\{ k_{\inmeshlg{o}} \text{ s.t. } k_{\inmeshlg{o}} \leftarrow j_{\inmeshlg{n}} \} \subseteq \{ k_{\inmeshlg{o}} \text{ s.t. } k_{\inmeshlg{o}} \rightarrow j_{\inmeshlg{n}} \}$, concluding the proof.

Secondly, we show that for a given $k_{\inmeshlg{o}}$
\nbea
\lvert \{ h_{\inmeshlg{n}} \text{ s.t. } k_{\inmeshlg{o}} \eitherarrow h_{\inmeshlg{n}} \}\rvert = 1.
\neea

\emph{Proof:}
We want to show that $$|\{ h_{\inmeshlg{n}} \text{ s.t. } k_{\inmeshlg{o}} \eitherarrow h_{\inmeshlg{n}} \}| = |\{ h_{\inmeshlg{n}} \text{ s.t. } k_{\inmeshlg{o}} \rightarrow h_{\inmeshlg{n}} \} \cup \{ h_{\inmeshlg{n}} \text{ s.t. } k_{\inmeshlg{o}} \leftarrow h_{\inmeshlg{n}} \}| = 1 .$$ If $\{ h_{\inmeshlg{n}} \text{ s.t. } k_{\inmeshlg{o}} \leftarrow h_{\inmeshlg{n}} \} = \emptyset$, this is trivially the case since $k_{\inmeshlg{o}}$ only has one nearest neighbour i.e. $\lvert \{ h_{\inmeshlg{n}} \text{ s.t. } k_{\inmeshlg{o}} \rightarrow h_{\inmeshlg{n}} \}\rvert = 1$. Therefore, in the following we consider only the case where it is not.

With this, we can show it is also the case that $|\{ h_{\inmeshlg{n}} \text{ s.t. } k_{\inmeshlg{o}} \leftarrow h_{\inmeshlg{n}} \}|=1$, which can be proven by contradiction. Suppose that $|\{ h_{\inmeshlg{n}} \text{ s.t. } k_{\inmeshlg{o}} \leftarrow h_{\inmeshlg{n}} \}|=R \neq 1$. This means that there exist multiple distinct nodes $\left\{h^r_{\inmeshlg{n}}\right\}_{r=1}^R$, each satisfying $k_{\inmeshlg{o}} \leftarrow h^r_{\inmeshlg{n}}$. Via the agglomerative condition in \autoref{eq:agg}, it follows that each of these nodes must also satisfy $k_{\inmeshlg{o}} \rightarrow h^r_{\inmeshlg{n}}$. This is not possible since $k_{\inmeshlg{o}}$ can only have a single nearest neighbour, meaning therefore $|\{ h_{\inmeshlg{n}} \text{ s.t. } k_{\inmeshlg{o}} \rightarrow h_{\inmeshlg{n}} \}|=1$.

We note that the single node in the set $|\{ h_{\inmeshlg{n}} \text{ s.t. } k_{\inmeshlg{o}} \leftarrow h_{\inmeshlg{n}} \}|$ must also satisfy $k_{\inmeshlg{o}} \rightarrow h_{\inmeshlg{n}}$ due to \autoref{eq:agg}. Therefore, it is in fact the exact same node as the node in the set $|\{ h_{\inmeshlg{n}} \text{ s.t. } k_{\inmeshlg{o}} \rightarrow h_{\inmeshlg{n}} \}|$, thus concluding the proof.

The updates can therefore be simplified to
\nbea
\tilde{W}^e_{ij_{\inmesh{n}}} &=& \sum_{\forall k_{\inmesh{o}} \text{ s.t. }  k_{\inmesh{o}} \rightarrow j_{\inmesh{n}}} W^e_{ik_{\inmesh{o}}}, \\
\tilde{W}^d_{i_{\inmesh{n}}j} &=& \underset{\forall k_{\inmesh{o}} \text{ s.t. }  k_{\inmesh{o}} \rightarrow i_{\inmesh{n}}}{\operatorname{mean}} W^d_{k_{\inmesh{o}}j}, \\
\tilde{b}^d_{i_{\inmesh{n}}} &=&  \underset{\forall k_{\inmesh{o}} \text{ s.t. }  k_{\inmesh{o}} \rightarrow i_{\inmesh{n}}}{\operatorname{mean}} b^d_{k_{\inmesh{o}}},
\neea
which therefore concludes the proof.

\subsection{Self-consistency}
\label{app:self_consistency}
We consider a sequence of two transforms, firstly an expansive transform $\mathcal{M}_o \to \mathcal{M}_n$, and then secondly a transform $\mathcal{M}_n \to \mathcal{M}_o$. Note that if $\mathcal{M}_o \to \mathcal{M}_n$ is expansive, then it must follow by definition that $\mathcal{M}_n \to \mathcal{M}_o$ is agglomerative.
The first transform $\mathcal{M}_o \to \mathcal{M}_n$ is expansive. Following \autoref{eq:exp}, this gives the first update as
\nbea
\hat{W}^e_{ij_{\inmesh{n}}} &=& \frac{1}{|\{h_{\inmeshlg{n}} \text{ s.t. } k_{\inmeshlg{o}} \leftarrow h_{\inmeshlg{n}} \}|} {W}^e_{ik_{\inmesh{o}}}, \qquad \text{ where}\quad k_{\inmeshlg{o}} \leftarrow j_{\inmeshlg{n}}, \\
\hat{W}^d_{i_{\inmesh{n}}j} &=& {W}^d_{k_{\inmesh{o}}j}, \qquad \text{ where}\quad k_{\inmeshlg{o}} \leftarrow i_{\inmeshlg{n}}, \\
\hat{b}^d_{i_{\inmesh{n}}} &=& {b}^d_{k_{\inmesh{o}}},  \qquad \text{ where}\quad k_{\inmeshlg{o}} \leftarrow i_{\inmeshlg{n}}.
\neea
The second transform $\mathcal{M}_o \leftarrow \mathcal{M}_n$ is agglomerative (note that the transform is \emph{not} $\mathcal{M}_o \to \mathcal{M}_n$). Following \autoref{eq:agg}, this gives the second update as
\nbea
\tilde{W}^e_{ij_{\inmesh{o}}} &=& \sum_{\forall k_{\inmesh{n}} \text{ s.t. } j_{\inmesh{o}} \leftarrow k_{\inmesh{n}}} \hat{W}^e_{ik_{\inmesh{n}}}, \\
\tilde{W}^d_{i_{\inmesh{o}}j} &=& \underset{\forall k_{\inmesh{n}} \text{ s.t. } i_{\inmesh{o}} \leftarrow k_{\inmesh{n}}}{\operatorname{mean}} \hat{W}^d_{k_{\inmesh{n}}j}, \\
\tilde{b}^d_{i_{\inmesh{o}}} &=&  \underset{\forall k_{\inmesh{n}} \text{ s.t. } i_{\inmesh{o}} \leftarrow k_{\inmesh{n}}}{\operatorname{mean}} b^d_{k_{\inmesh{n}}}.
\neea
Subbing everything in gives
\nbea
\tilde{W}^e_{ij_{\inmesh{o}}} &=& \sum_{\forall k_{\inmesh{n}} \text{ s.t. } j_{\inmesh{o}} \leftarrow k_{\inmesh{n}}} \frac{1}{|\{h_{\inmeshlg{n}} \text{ s.t. } l_{\inmeshlg{o}} \leftarrow h_{\inmeshlg{n}} \}|} {W}^e_{il_{\inmesh{o}}}, \qquad \text{ where}\quad l_{\inmeshlg{o}} \leftarrow k_{\inmeshlg{n}},\\
&=& \sum_{\forall k_{\inmesh{n}} \text{ s.t. } j_{\inmesh{o}} \leftarrow k_{\inmesh{n}}} \frac{1}{|\{h_{\inmeshlg{n}} \text{ s.t. } j_{\inmeshlg{o}} \leftarrow h_{\inmeshlg{n}} \}|} {W}^e_{ij_{\inmesh{o}}},\\
&=& W^e_{ij_{\inmesh{o}}},\\
\tilde{W}^d_{i_{\inmesh{o}}j} &=& \underset{\forall k_{\inmesh{n}} \text{ s.t. } i_{\inmesh{o}} \leftarrow k_{\inmesh{n}}}{\operatorname{mean}} {W}^d_{l_{\inmesh{o}}j}, \qquad \text{ where}\quad l_{\inmeshlg{o}} \leftarrow k_{\inmeshlg{n}}, \\
&=& \underset{\forall k_{\inmesh{n}} \text{ s.t. } i_{\inmesh{o}} \leftarrow k_{\inmesh{n}}}{\operatorname{mean}} {W}^d_{i_{\inmesh{o}}j},\\
&=& W^d_{i_{\inmesh{o}}j},\\
\tilde{b}^d_{i_{\inmesh{o}}} &=&  \underset{\forall k_{\inmesh{n}} \text{ s.t. } i_{\inmesh{o}} \leftarrow k_{\inmesh{n}}}{\operatorname{mean}} {b}^d_{l_{\inmesh{o}}},  \qquad \text{ where}\quad l_{\inmeshlg{o}} \leftarrow k_{\inmeshlg{n}},\\
&=& \underset{\forall k_{\inmesh{n}} \text{ s.t. } i_{\inmesh{o}} \leftarrow k_{\inmesh{n}}}{\operatorname{mean}} {b}^d_{i_{\inmesh{o}}},\\
&=& b^d_{i_{\inmesh{o}}}.
\neea
This concludes the proof.

\subsection{General GFN transform as a composition of an expansive and agglomerative transform}
\label{app:gfn_decomp}
Splitting up the summations, the general GFN transforms given in \autoref{eq:gfn} can be rewritten as
\nbea
\begin{aligned}
\tilde{W}^e_{ij_{\inmesh{n}}} &= \left[ \sum_{\forall k_{\inmesh{o}} \text{ s.t } k_{\inmesh{o}} \rightarrow j_{\inmesh{n}}}
+ \sum\limits_{\substack{\forall k_{\inmesh{o}} \text{ s.t } k_{\inmesh{o}} \leftarrow j_{\inmesh{n}} \\ \text{but not}~k_{\inmesh{o}} \rightarrow j_{\inmesh{n}}}} \right] \frac{W^e_{ik_{\inmesh{o}}}}{\lvert \{ h_{\inmeshlg{n}} \text{ s.t. } k_{\inmeshlg{o}} ~\eitherarrow~ h_{\inmeshlg{n}} \}\rvert}, \\
\tilde{W}^d_{i_{\inmesh{n}}j} &= \frac{1}{|\{ k_{\inmeshlg{o}} \text{ s.t } k_{\inmeshlg{o}} \eitherarrow i_{\inmeshlg{n}} \}|} \left[ \sum_{\forall k_{\inmesh{o}} \text{ s.t } k_{\inmesh{o}} \rightarrow i_{\inmesh{n}}} + \sum\limits_{\substack{\forall k_{\inmesh{o}} \text{ s.t } k_{\inmesh{o}} \leftarrow i_{\inmesh{n}} \\ \text{but not}~k_{\inmesh{o}} \rightarrow i_{\inmesh{n}}}} \right] {W}^d_{k_{\inmesh{o}}j}, \\
\tilde{b}^d_{i_{\inmesh{n}}} &= \frac{1}{|\{ k_{\inmeshlg{o}} \text{ s.t } k_{\inmeshlg{o}} \eitherarrow i_{\inmeshlg{n}} \}|} \left[ \sum_{\forall k_{\inmesh{o}} \text{ s.t } k_{\inmesh{o}} \rightarrow i_{\inmesh{n}}} + \sum\limits_{\substack{\forall k_{\inmesh{o}} \text{ s.t } k_{\inmesh{o}} \leftarrow i_{\inmesh{n}} \\ \text{but not}~k_{\inmesh{o}} \rightarrow i_{\inmesh{n}}}} \right] b^d_{k_{\inmesh{o}}}.
\end{aligned}
\neea
By defining an auxiliary mesh
\nbea
\mathcal{M}_{} = \mathcal{M}_o \cup \{\mathcal{M}_n \left[i_{\inmeshlg{n}}\right] \text{ s.t. } j_{\inmeshlg{o}} \leftarrow i_{\inmeshlg{n}} \text{ but not } j_{\inmeshlg{o}} \rightarrow i_{\inmeshlg{n}} \},
\neea
we can convert the summations to
\nbea
\sum_{\forall k_{\inmesh{o}} \text{ s.t } k_{\inmesh{o}} \rightarrow i_{\inmesh{n}}} \longrightarrow \sum\limits_{\substack{\forall l_{\inmesh{}} \text{ s.t } l_{\inmesh{}} \rightarrow i_{\inmesh{n}} \\ \text{and}~l_{\inmesh{}} \in \mathcal{M}_o}}, \quad k_{\inmeshlg{o}} \leftarrow l_{\inmeshlg{}}, \\
\sum\limits_{\substack{\forall k_{\inmesh{o}} \text{ s.t } k_{\inmesh{o}} \leftarrow i_{\inmesh{n}} \\ \text{but not}~k_{\inmesh{o}} \rightarrow i_{\inmesh{n}}}} \longrightarrow \sum\limits_{\substack{\forall l_{\inmesh{}} \text{ s.t } l_{\inmesh{}} \rightarrow i_{\inmesh{n}} \\ \text{and}~l_{\inmesh{}} \in \mathcal{M}_n}}, \quad k_{\inmeshlg{o}} \leftarrow l_{\inmeshlg{}}.
\neea
Subbing this in, this gives an updated expression for the general GFN transforms of
\nbea
\begin{aligned}
\tilde{W}^e_{ij_{\inmesh{n}}} &= \left[ \sum\limits_{\substack{\forall l_{\inmesh{}} \text{ s.t } l_{\inmesh{}} \rightarrow j_{\inmesh{n}} \\ \text{and}~l_{\inmesh{}} \in \mathcal{M}_o}}
+ \sum\limits_{\substack{\forall l_{\inmesh{}} \text{ s.t } l_{\inmesh{}} \rightarrow j_{\inmesh{n}} \\ \text{and}~l_{\inmesh{}} \in \mathcal{M}_n}} \right] \frac{W^e_{ik_{\inmesh{o}}}}{\lvert \{ h_{\inmeshlg{n}} \text{ s.t. } k_{\inmeshlg{o}} ~\eitherarrow~ h_{\inmeshlg{n}} \}\rvert}, \quad  k_{\inmeshlg{o}} \leftarrow l_{\inmeshlg{}},\\
\tilde{W}^d_{i_{\inmesh{n}}j} &= \frac{1}{|\{ g_{\inmeshlg{o}} \text{ s.t } g_{\inmeshlg{o}} \eitherarrow i_{\inmeshlg{n}} \}|} \left[ \sum\limits_{\substack{\forall l_{\inmesh{}} \text{ s.t } l_{\inmesh{}} \rightarrow i_{\inmesh{n}} \\ \text{and}~l_{\inmesh{}} \in \mathcal{M}_o}}
+ \sum\limits_{\substack{\forall l_{\inmesh{}} \text{ s.t } l_{\inmesh{}} \rightarrow i_{\inmesh{n}} \\ \text{and}~l_{\inmesh{}} \in \mathcal{M}_n}} \right] {W}^d_{k_{\inmesh{o}}j}, \quad  k_{\inmeshlg{o}} \leftarrow l_{\inmeshlg{}},\\
\tilde{b}^d_{i_{\inmesh{n}}} &= \frac{1}{|\{ g_{\inmeshlg{o}} \text{ s.t } g_{\inmeshlg{o}} \eitherarrow i_{\inmeshlg{n}} \}|} \left[ \sum\limits_{\substack{\forall l_{\inmesh{}} \text{ s.t } l_{\inmesh{}} \rightarrow i_{\inmesh{n}} \\ \text{and}~l_{\inmesh{}} \in \mathcal{M}_o}}
+ \sum\limits_{\substack{\forall l_{\inmesh{m}} \text{ s.t } l_{\inmesh{}} \rightarrow i_{\inmesh{n}} \\ \text{and}~l_{\inmesh{}} \in \mathcal{M}_n}} \right] {b}^d_{i_{\inmesh{n}}}, \quad  k_{\inmeshlg{o}} \leftarrow l_{\inmeshlg{}}.
\end{aligned}
\neea
Since it is the case by construction of $\mathcal{M}$ that $|\{ h_{\inmeshlg{n}} \text{ s.t. } k_{\inmeshlg{o}} ~\eitherarrow~ h_{\inmeshlg{n}} \}| = |\{ h_{\inmeshlg{}} \text{ s.t. } k_{\inmeshlg{o}} \leftarrow h_{\inmeshlg{}} \}|$ and $|\{ g_{\inmeshlg{o}} \text{ s.t } g_{\inmeshlg{o}} \eitherarrow i_{\inmeshlg{n}} \}| = |\{ l_{\inmeshlg{}} \text{ s.t } l_{\inmeshlg{}} \rightarrow i_{\inmeshlg{n}} \}|$, we can further rewrite these expressions as
\nbea
\begin{aligned}
\tilde{W}^e_{ij_{\inmesh{n}}} &= \sum_{\substack{\forall l_{\inmesh{}} \text{ s.t } l_{\inmesh{}} \rightarrow j_{\inmesh{n}}}} \frac{W^e_{ik_{\inmesh{o}}}}{\lvert \{ h_{\inmeshlg{}} \text{ s.t. } k_{\inmeshlg{o}} \leftarrow h_{\inmeshlg{}} \}\rvert}, \quad  k_{\inmeshlg{o}} \leftarrow l_{\inmeshlg{}},\\
\tilde{W}^d_{i_{\inmesh{n}}j} &= \underset{\forall l_{\inmesh{}} \text{ s.t } l_{\inmesh{}} \rightarrow i_{\inmesh{n}}}{\operatorname{mean}} {W}^d_{k_{\inmesh{o}}j}, \quad  k_{\inmeshlg{o}} \leftarrow l_{\inmeshlg{}},\\
\tilde{b}^d_{i_{\inmesh{n}}} &= \underset{\forall l_{\inmesh{}} \text{ s.t } l_{\inmesh{}} \rightarrow i_{\inmesh{n}}}{\operatorname{mean}} {b}^d_{k_{\inmesh{o}}}, \quad  k_{\inmeshlg{o}} \leftarrow l_{\inmeshlg{}}.
\end{aligned}
\neea
We can finally compute the transform in two steps, as
\begin{equation*}
\begin{rcases}
\hat{W}^e_{ij_{\inmesh{}}} &= \frac{1}{|\{h_{\inmesh{}} \text{ s.t. } k_{\inmesh{o}} \leftarrow h_{\inmesh{}} \}|} {W}^e_{ik_{\inmesh{o}}}, \qquad \text{ where}\quad k_{\inmeshlg{o}} \leftarrow j_{\inmeshlg{}}, \\
\hat{W}^d_{i_{\inmesh{}}j} &= {W}^d_{k_{\inmesh{o}}j}, \qquad \text{ where}\quad k_{\inmeshlg{o}} \leftarrow i_{\inmeshlg{}}, \\
{{\hat{b}}^d}_{i_{\inmesh{}}} &= {b}^d_{k_{\inmesh{o}}},  \qquad \text{ where}\quad k_{\inmeshlg{o}} \leftarrow i_{\inmeshlg{}},
\end{rcases} \quad \text{Expansion step}
\end{equation*}
followed by
\begin{equation*}
\begin{rcases}
\tilde{W}^e_{ij_{\inmesh{n}}} &= \sum_{\forall k_{\inmesh{}} \text{ s.t. }  k_{\inmesh{}} \rightarrow j_{\inmesh{n}}} {\hat{W}}^e_{ik_{\inmesh{}}}, \\
\tilde{W}^d_{i_{\inmesh{n}}j} &= \underset{\forall k_{\inmesh{}} \text{ s.t. }  k_{\inmesh{}} \rightarrow i_{\inmesh{n}}}{\operatorname{mean}} {\hat{W}}^d_{k_{\inmesh{}}j}, \\
\tilde{b}^d_{i_{\inmesh{n}}} &=  \underset{\forall k_{\inmesh{}} \text{ s.t. }  k_{\inmesh{}} \rightarrow i_{\inmesh{n}}}{\operatorname{mean}} {\hat{b}}^d_{k_{\inmesh{}}},
\end{rcases} \quad \text{Agglomerative step}
\end{equation*}
which concludes the proof.

\subsection{Super- and sub-resolution error bound for GFN-ROM}
\label{subsec:error_gfnrom}
We seek a bound on the following quantity
\nbea
\left| u(x_{i_{\inmesh{n}}}) - \operatorname{dec}^{\mathcal{M}_o \to \mathcal{M}_n}(\operatorname{map}(\boldsymbol{\mu}))_{i_{\inmesh{n}}} \right|, \forall i_{\inmeshlg{n}}.
\neea
We can rewrite the expression and bound it as
\nbea
&& \left| u(x_{i_{\inmesh{n}}}) - \operatorname{dec}^{\mathcal{M}_o \to \mathcal{M}_n}(\operatorname{map}(\boldsymbol{\mu}))_{i_{\inmesh{n}}} \right| \\
&=& \left| u(x_{i_{\inmesh{n}}}) - \sum_{j=1}^{L} \tilde{W}^d_{i_{\inmesh{n}}j} \operatorname{map}(\boldsymbol{\mu})_j + \tilde{b}^d_{i_{\inmesh{n}}} \right|,\\
&=& \left| u(x_{i_{\inmesh{n}}}) - \sum_{j=1}^{L} \underset{\forall k_{\inmesh{o}} \text{ s.t } k_{\inmesh{o}} \eitherarrow i_{\inmesh{n}}}{\operatorname{mean}} {W}^d_{k_{\inmesh{o}}j} \operatorname{map}(\boldsymbol{\mu})_j + \underset{\forall k_{\inmesh{o}} \text{ s.t } k_{\inmesh{o}} \eitherarrow i_{\inmesh{n}}}{\operatorname{mean}} {b}^d_{k_{\inmesh{o}}} \right|, \\
&=& \left| u(x_{i_{\inmesh{n}}}) - \underset{\forall k_{\inmesh{o}} \text{ s.t } k_{\inmesh{o}} \eitherarrow i_{\inmesh{n}}}{\operatorname{mean}} \left( \sum_{j=1}^{L} {W}^d_{k_{\inmesh{o}}j} \operatorname{map}(\boldsymbol{\mu})_j + {b}^d_{k_{\inmesh{o}}} \right) \right| \\
&=& \left| u(x_{i_{\inmesh{n}}}) - \underset{\forall k_{\inmesh{o}} \text{ s.t } k_{\inmesh{o}} \eitherarrow i_{\inmesh{n}}}{\operatorname{mean}} \operatorname{dec}^{\mathcal{M}_o \to \mathcal{M}_o}(\operatorname{map}(\boldsymbol{\mu}))_{k_{\inmesh{o}}} \right|, \\
&=& \left| \underset{\forall k_{\inmesh{o}} \text{ s.t } k_{\inmesh{o}} \eitherarrow i_{\inmesh{n}}}{\operatorname{mean}} \left( u(x_{i_{\inmesh{n}}}) - u(x_{k_{\inmesh{o}}}) + u(x_{k_{\inmesh{o}}}) - \operatorname{dec}^{\mathcal{M}_o \to \mathcal{M}_o}(\operatorname{map}(\boldsymbol{\mu}))_{k_{\inmesh{o}}} \right) \right|, \\
&\leq& \underset{\forall k_{\inmesh{o}} \text{ s.t } k_{\inmesh{o}} \eitherarrow i_{\inmesh{n}}}{\operatorname{mean}} \left| u(x_{i_{\inmesh{n}}}) - u(x_{k_{\inmesh{o}}}) + u(x_{k_{\inmesh{o}}}) - \operatorname{dec}^{\mathcal{M}_o \to \mathcal{M}_o}(\operatorname{map}(\boldsymbol{\mu}))_{k_{\inmesh{o}}} \right|, \\
&\leq& \underset{\forall k_{\inmesh{o}} \text{ s.t } k_{\inmesh{o}} \eitherarrow i_{\inmesh{n}}}{\operatorname{mean}} \left( \left| u(x_{i_{\inmesh{n}}}) - u(x_{k_{\inmesh{o}}}) \right| + \left| u(x_{k_{\inmesh{o}}}) - \operatorname{dec}^{\mathcal{M}_o \to \mathcal{M}_o}(\operatorname{map}(\boldsymbol{\mu}))_{k_{\inmesh{o}}} \right| \right), \\
&\leq& \underset{\forall k_{\inmesh{o}} \text{ s.t } k_{\inmesh{o}} \eitherarrow i_{\inmesh{n}}}{\operatorname{mean}} \left( \delta + \tau \right),\\
&=& \delta + \tau.
\neea
This concludes the proof. For the multiple layer case, the result is also the same.

\subsection{Super- and sub-resolution error bound for GFN-ROM mapper}
\label{subsec:error_mapper}
We seek a bound on
\bea
\left| \operatorname{map}(\boldsymbol{\mu})_i - \operatorname{enc}^{\mathcal{M}_o \to \mathcal{M}_n}(\bfmu_{\inmeshlg{n}})_i \right|, \forall i.
\eea
We can rewrite the expression as
\nbea
&& \left| \operatorname{map}(\boldsymbol{\mu})_i - \operatorname{enc}^{\mathcal{M}_o \to \mathcal{M}_o}(\bfmu_{\inmeshlg{o}})_i + \operatorname{enc}^{\mathcal{M}_o \to \mathcal{M}_o}(\bfmu_{\inmeshlg{o}})_i - \operatorname{enc}^{\mathcal{M}_o \to \mathcal{M}_n}(\bfmu_{\inmeshlg{n}})_i \right| \\
&\leq& \left| \operatorname{map}(\boldsymbol{\mu})_i - \operatorname{enc}^{\mathcal{M}_o \to \mathcal{M}_o}(\bfmu_{\inmeshlg{o}})_i \right| + \left| \operatorname{enc}^{\mathcal{M}_o \to \mathcal{M}_o}(\bfmu_{\inmeshlg{o}})_i - \operatorname{enc}^{\mathcal{M}_o \to \mathcal{M}_n}(\bfmu_{\inmeshlg{n}})_i \right|,\\
&\leq& \alpha + \left| \operatorname{enc}^{\mathcal{M}_o \to \mathcal{M}_o}(\bfmu_{\inmeshlg{o}})_i - \operatorname{enc}^{\mathcal{M}_o \to \mathcal{M}_n}(\bfmu_{\inmeshlg{n}})_i \right|,\\
&=& \alpha + \left| \sigma \left( \sum_{j_{\inmesh{o}}=1}^{N_o} W_{ij_{\inmesh{o}}}^e u(x_{j_{\inmesh{o}}}) + b^e_{i} \right) - \sigma \left(\sum_{k_{\inmesh{n}}=1}^{N_n} \tilde{W}_{k_{\inmesh{n}}i}^e u(x_{k_{\inmesh{n}}}) + {\tilde{b}^e}_{i}\right) \right|.
\neea
We assume a Lipschitz condition with Lipschitz constant $C$ for the activation function $\sigma$, which allows us to continue bounding as
\nbea
&\leq& \alpha + C \left| \sum_{j_{\inmesh{o}}=1}^{N_o} W_{ij_{\inmesh{o}}}^e u(x_{j_{\inmesh{o}}}) + b^e_{i} - \left(\sum_{k_{\inmesh{n}}=1}^{N_n} \tilde{W}_{ik_{\inmesh{n}}}^e u(x_{k_{\inmesh{n}}}) + {\tilde{b}^e}_{i}\right) \right|, \\
&=& \alpha + C \left| \sum_{j_{\inmesh{o}}=1}^{N_o} W_{ij_{\inmesh{o}}}^e u(x_{j_{\inmesh{o}}}) - \sum_{k_{\inmesh{n}}=1}^{N_n} \tilde{W}_{ik_{\inmesh{n}}}^e u(x_{k_{\inmesh{n}}}) \right|, \\
&=& \alpha + C \left| \sum_{j_{\inmesh{o}}=1}^{N_o} W_{ij_{\inmesh{o}}}^e u(x_{j_{\inmesh{o}}}) - \sum_{k_{\inmesh{n}}=1}^{N_n} ~ \sum_{\forall m_{\inmesh{o}} \text{ s.t } m_{\inmesh{o}} \eitherarrow k_{\inmesh{n}}} \frac{W^e_{im_{\inmesh{o}}}}{\lvert \{ h_{\inmeshlg{n}} \text{ s.t. } m_{\inmeshlg{o}} ~\eitherarrow~ h_{\inmeshlg{n}} \}\rvert} u(x_{k_{\inmesh{n}}}) \right|, \\
&=& \alpha + C \left| \sum_{j_{\inmesh{o}}=1}^{N_o} W_{ij_{\inmesh{o}}}^e \left( u(x_{j_{\inmesh{o}}}) - \sum_{\forall k_{\inmesh{n}} \text{ s.t } j_{\inmesh{o}} \eitherarrow k_{\inmesh{n}}} \frac{u(x_{k_{\inmesh{n}}})}{\lvert \{ h_{\inmeshlg{n}} \text{ s.t. } j_{\inmeshlg{o}} ~\eitherarrow~ h_{\inmeshlg{n}} \}\rvert} \right) \right|, \\
&=& \alpha + C \left| \sum_{j_{\inmesh{o}}=1}^{N_o} \frac{W_{ij_{\inmesh{o}}}^e}{\lvert \{ h_{\inmeshlg{n}} \text{ s.t. } j_{\inmeshlg{o}} ~\eitherarrow~ h_{\inmeshlg{n}} \}\rvert} \sum_{\forall k_{\inmesh{n}} \text{ s.t } j_{\inmesh{o}} \eitherarrow k_{\inmesh{n}}} \left(  u(x_{j_{\inmesh{o}}}) - u(x_{k_{\inmesh{n}}}) \right) \right|, \\
&=& \alpha + C \left| \sum_{j_{\inmesh{o}}=1}^{N_o} \frac{W_{ij_{\inmesh{o}}}^e}{\lvert \{ h_{\inmeshlg{n}} \text{ s.t. } j_{\inmeshlg{o}} ~\eitherarrow~ h_{\inmeshlg{n}} \}\rvert} \sum_{\forall k_{\inmesh{n}} \text{ s.t } j_{\inmesh{o}} \eitherarrow k_{\inmesh{n}}} \delta \right|, \\
&\leq& \alpha + \delta C \sum_{j_{\inmesh{o}}=1}^{N_o} |W_{ij_{\inmesh{o}}}^e|.
\neea
To obtain a bound independent of $i$, one can bound as
\nbea
\leq \alpha + \delta C ||\bfmW^{e}||_\infty.
\neea
This concludes the proof. For the multiple layer case, the result is
\nbea
&\leq& \alpha + \delta C^{P+1} \sum_{l_{P}=1}^{L_{P}} |W^{(P)}_{i,l_{P}}|  \sum_{l_{P-1}=1}^{L_{P-1}} |W^{(P-1)}_{l_{P},l_{P-1}}| \cdots \sum_{l_{1}=1}^{L_{1}} |W^{(1)}_{l_{2},l_{1}}| ~ \sum_{j_{\inmesh{o}}=1}^{N_o} |W_{l_1,j_{\inmesh{o}}}^e|,\\
&\leq& \alpha + \delta C^{P+1} ||\bfmW^{e}||_\infty \prod_{p=1}^P ||\bfmW^{(p)}||_\infty.
\neea

\subsection{Super- and sub-resolution error bound for GFN-ROM autoencoder}
\label{subsec:error_autoenc}
We seek a bound on
\nbea
\left| u(x_{i_{\inmesh{n}}}) - \operatorname{dec}^{\mathcal{M}_o \to \mathcal{M}_n}(\operatorname{enc}^{\mathcal{M}_o \to \mathcal{M}_n}(\bfmu_{\inmeshlg{n}}))_{i_{\inmesh{n}}} \right|, \forall i_{\inmeshlg{n}}.
\neea
To begin bounding this, we can rewrite the expression as
\nbea
&=& \left| u(x_{i_{\inmesh{n}}}) - \underset{\forall k_{\inmesh{o}} \text{ s.t } k_{\inmesh{o}} \eitherarrow i_{\inmesh{n}}}{\operatorname{mean}} u(x_{k_{\inmesh{o}}}) \right. \\
&& + \underset{\forall k_{\inmesh{o}} \text{ s.t } k_{\inmesh{o}} \eitherarrow i_{\inmesh{n}}}{\operatorname{mean}} u(x_{k_{\inmesh{o}}}) - \underset{\forall k_{\inmesh{o}} \text{ s.t } k_{\inmesh{o}} \eitherarrow i_{\inmesh{n}}}{\operatorname{mean}} \operatorname{dec}^{\mathcal{M}_o \to \mathcal{M}_o}(\operatorname{enc}^{\mathcal{M}_o \to \mathcal{M}_o}(\bfmu_{\inmeshlg{o}}))_{k_{\inmesh{o}}} \\
&& \left. + \underset{\forall k_{\inmesh{o}} \text{ s.t } k_{\inmesh{o}} \eitherarrow i_{\inmesh{n}}}{\operatorname{mean}} \operatorname{dec}^{\mathcal{M}_o \to \mathcal{M}_o}(\operatorname{enc}^{\mathcal{M}_o \to \mathcal{M}_o}(\bfmu_{\inmeshlg{o}}))_{k_{\inmesh{o}}} - \operatorname{dec}^{\mathcal{M}_o \to \mathcal{M}_n}(\operatorname{enc}^{\mathcal{M}_o \to \mathcal{M}_n}(\bfmu_{\inmeshlg{n}}))_{i_{\inmesh{n}}} \right|, \\
&\leq& \left| u(x_{i_{\inmesh{n}}}) - \underset{\forall k_{\inmesh{o}} \text{ s.t } k_{\inmesh{o}} \eitherarrow i_{\inmesh{n}}}{\operatorname{mean}} u(x_{k_{\inmesh{o}}}) \right| \\
&& + \left| \underset{\forall k_{\inmesh{o}} \text{ s.t } k_{\inmesh{o}} \eitherarrow i_{\inmesh{n}}}{\operatorname{mean}} u(x_{k_{\inmesh{o}}}) - \underset{\forall k_{\inmesh{o}} \text{ s.t } k_{\inmesh{o}} \eitherarrow i_{\inmesh{n}}}{\operatorname{mean}} \operatorname{dec}^{\mathcal{M}_o \to \mathcal{M}_o}(\operatorname{enc}^{\mathcal{M}_o \to \mathcal{M}_o}(\bfmu_{\inmeshlg{o}}))_{k_{\inmesh{o}}} \right| \\
&& + \left| \underset{\forall k_{\inmesh{o}} \text{ s.t } k_{\inmesh{o}} \eitherarrow i_{\inmesh{n}}}{\operatorname{mean}} \operatorname{dec}^{\mathcal{M}_o \to \mathcal{M}_o}(\operatorname{enc}^{\mathcal{M}_o \to \mathcal{M}_o}(\bfmu_{\inmeshlg{o}}))_{k_{\inmesh{o}}} - \operatorname{dec}^{\mathcal{M}_o \to \mathcal{M}_n}(\operatorname{enc}^{\mathcal{M}_o \to \mathcal{M}_n}(\bfmu_{\inmeshlg{n}}))_{i_{\inmesh{n}}} \right|, \\
&\leq& \delta + \beta + \left| \underset{\forall k_{\inmesh{o}} \text{ s.t } k_{\inmesh{o}} \eitherarrow i_{\inmesh{n}}}{\operatorname{mean}} \operatorname{dec}^{\mathcal{M}_o \to \mathcal{M}_o}(\operatorname{enc}^{\mathcal{M}_o \to \mathcal{M}_o}(\bfmu_{\inmeshlg{o}}))_{k_{\inmesh{o}}} - \operatorname{dec}^{\mathcal{M}_o \to \mathcal{M}_n}(\operatorname{enc}^{\mathcal{M}_o \to \mathcal{M}_n}(\bfmu_{\inmeshlg{n}}))_{i_{\inmesh{n}}} \right|.
\neea
We consider the last term in more detail. We have
\nbea
&& \left| \underset{\forall k_{\inmesh{o}} \text{ s.t } k_{\inmesh{o}} \eitherarrow i_{\inmesh{n}}}{\operatorname{mean}} \operatorname{dec}^{\mathcal{M}_o \to \mathcal{M}_o}(\operatorname{enc}^{\mathcal{M}_o \to \mathcal{M}_o}(\bfmu_{\inmeshlg{o}}))_{k_{\inmesh{o}}} - \operatorname{dec}^{\mathcal{M}_o \to \mathcal{M}_n}(\operatorname{enc}^{\mathcal{M}_o \to \mathcal{M}_n}(\bfmu_{\inmeshlg{n}}))_{i_{\inmesh{n}}} \right| \\
&=& \left| \underset{\forall k_{\inmesh{o}} \text{ s.t } k_{\inmesh{o}} \eitherarrow i_{\inmesh{n}}}{\operatorname{mean}} \left( \sum_{j=1}^{L} W^d_{k_{\inmesh{o}}j} \operatorname{enc}^{\mathcal{M}_o \to \mathcal{M}_o}(\bfmu_{\inmeshlg{o}})_{j} + b^d_{k_{\inmesh{o}}} \right) \right.\\
&& \left. - \sum_{m=1}^{L} \tilde{W}^d_{i_{\inmesh{n}}m} \operatorname{enc}^{\mathcal{M}_o \to \mathcal{M}_n}(\bfmu_{\inmeshlg{n}})_{m} - \tilde{b}^d_{i_{\inmesh{n}}} \right|, \\
&=& \left| \underset{\forall k_{\inmesh{o}} \text{ s.t } k_{\inmesh{o}} \eitherarrow i_{\inmesh{n}}}{\operatorname{mean}} \left( \sum_{j=1}^{L} W^d_{k_{\inmesh{o}}j} \operatorname{enc}^{\mathcal{M}_o \to \mathcal{M}_o}(\bfmu_{\inmeshlg{o}})_{j} + b^d_{k_{\inmesh{o}}} \right) \right.\\
&& \left. - \sum_{m=1}^{L} \underset{\forall k_{\inmesh{o}} \text{ s.t } k_{\inmesh{o}} \eitherarrow i_{\inmesh{n}}}{\operatorname{mean}} \left( W^d_{k_{\inmesh{o}}m} \right) \operatorname{enc}^{\mathcal{M}_o \to \mathcal{M}_n}(\bfmu_{\inmeshlg{n}})_{m} - \underset{\forall k_{\inmesh{o}} \text{ s.t } k_{\inmesh{o}} \eitherarrow i_{\inmesh{n}}}{\operatorname{mean}} \left( b^d_{k_{\inmesh{o}}} \right) \right|, \\ 
&=& \left| \underset{\forall k_{\inmesh{o}} \text{ s.t } k_{\inmesh{o}} \eitherarrow i_{\inmesh{n}}}{\operatorname{mean}} \sum_{j=1}^{L} W^d_{k_{\inmesh{o}}j} \left( \operatorname{enc}^{\mathcal{M}_o \to \mathcal{M}_o}(\bfmu_{\inmeshlg{o}})_{j} -\operatorname{enc}^{\mathcal{M}_o \to \mathcal{M}_n}(\bfmu_{\inmeshlg{n}})_{j}\right) \right|,\\
&\leq& \underset{\forall k_{\inmesh{o}} \text{ s.t } k_{\inmesh{o}} \eitherarrow i_{\inmesh{n}}}{\operatorname{mean}} \sum_{j=1}^{L} \left| W^d_{k_{\inmesh{o}}j} \right| \left| \operatorname{enc}^{\mathcal{M}_o \to \mathcal{M}_o}(\bfmu_{\inmeshlg{o}})_{j} -\operatorname{enc}^{\mathcal{M}_o \to \mathcal{M}_n}(\bfmu_{\inmeshlg{n}})_{j} \right|.
\neea
We assume a Lipschitz condition with Lipschitz constant $C$ for the activation function $\sigma$, allowing us to bound
\nbea
&& \left| \operatorname{enc}^{\mathcal{M}_o \to \mathcal{M}_o}(\bfmu_{\inmeshlg{o}})_{i} -\operatorname{enc}^{\mathcal{M}_o \to \mathcal{M}_n}(\bfmu_{\inmeshlg{n}})_{i} \right|\\
&=& \left| \sigma \left( \sum_{j_{\inmesh{o}}=1}^{N_o} W_{ij_{\inmesh{o}}}^e u(x_{j_{\inmesh{o}}}) + b^e_{i} \right) - \sigma \left(\sum_{k_{\inmesh{n}}=1}^{N_n} \tilde{W}_{k_{\inmesh{n}}i}^e u(x_{k_{\inmesh{n}}}) + {\tilde{b}^e}_{i}\right) \right|,\\
&\leq& C \left| \sum_{j_{\inmesh{o}}=1}^{N_o} W_{ij_{\inmesh{o}}}^e u(x_{j_{\inmesh{o}}}) + b^e_{i} - \left(\sum_{k_{\inmesh{n}}=1}^{N_n} \tilde{W}_{ik_{\inmesh{n}}}^e u(x_{k_{\inmesh{n}}}) + {\tilde{b}^e}_{i}\right) \right|, \\
&=& C \left| \sum_{j_{\inmesh{o}}=1}^{N_o} W_{ij_{\inmesh{o}}}^e u(x_{j_{\inmesh{o}}}) - \sum_{k_{\inmesh{n}}=1}^{N_n} \tilde{W}_{ik_{\inmesh{n}}}^e u(x_{k_{\inmesh{n}}}) \right|, \\
&=& C \left| \sum_{j_{\inmesh{o}}=1}^{N_o} W_{ij_{\inmesh{o}}}^e u(x_{j_{\inmesh{o}}}) - \sum_{k_{\inmesh{n}}=1}^{N_n} ~ \sum_{\forall m_{\inmesh{o}} \text{ s.t } m_{\inmesh{o}} \eitherarrow k_{\inmesh{n}}} \frac{W^e_{im_{\inmesh{o}}}}{\lvert \{ h_{\inmeshlg{n}} \text{ s.t. } m_{\inmeshlg{o}} ~\eitherarrow~ h_{\inmeshlg{n}} \}\rvert} u(x_{k_{\inmesh{n}}}) \right|, \\
&=& C \left| \sum_{j_{\inmesh{o}}=1}^{N_o} W_{ij_{\inmesh{o}}}^e \left( u(x_{j_{\inmesh{o}}}) - \sum_{\forall k_{\inmesh{n}} \text{ s.t } j_{\inmesh{o}} \eitherarrow k_{\inmesh{n}}} \frac{u(x_{k_{\inmesh{n}}})}{\lvert \{ h_{\inmeshlg{n}} \text{ s.t. } j_{\inmeshlg{o}} ~\eitherarrow~ h_{\inmeshlg{n}} \}\rvert} \right) \right|, \\
&=& C \left| \sum_{j_{\inmesh{o}}=1}^{N_o} \frac{W_{ij_{\inmesh{o}}}^e}{\lvert \{ h_{\inmeshlg{n}} \text{ s.t. } j_{\inmeshlg{o}} ~\eitherarrow~ h_{\inmeshlg{n}} \}\rvert} \sum_{\forall k_{\inmesh{n}} \text{ s.t } j_{\inmesh{o}} \eitherarrow k_{\inmesh{n}}} \left(  u(x_{j_{\inmesh{o}}}) - u(x_{k_{\inmesh{n}}}) \right) \right|, \\
&=& C \left| \sum_{j_{\inmesh{o}}=1}^{N_o} \frac{W_{ij_{\inmesh{o}}}^e}{\lvert \{ h_{\inmeshlg{n}} \text{ s.t. } j_{\inmeshlg{o}} ~\eitherarrow~ h_{\inmeshlg{n}} \}\rvert} \sum_{\forall k_{\inmesh{n}} \text{ s.t } j_{\inmesh{o}} \eitherarrow k_{\inmesh{n}}} \delta \right|, \\
&\leq& \delta C \sum_{j_{\inmesh{o}}=1}^{N_o} |W_{ij_{\inmesh{o}}}^e|.
\neea
Substituting these expressions into the original bound, we obtain
\nbea
&\leq& \delta + \beta + \delta C \underset{\forall k_{\inmesh{o}} \text{ s.t } k_{\inmesh{o}} \eitherarrow i_{\inmesh{n}}}{\operatorname{mean}} \sum_{j=1}^{L} \left| W^d_{k_{\inmesh{o}}j} \right| \sum_{m_{\inmesh{o}}=1}^{N_o} \left|W_{jm_{\inmesh{o}}}^e\right|.
\neea
To obtain a bound independent of $i$, one can bound as
\nbea
&\leq& \delta + \beta + \delta C ||\bfmW^{d}||_\infty ||\bfmW^{e}||_\infty.
\neea
This concludes the proof. For the multiple layer case, the result is
\nbea
&\leq& \beta + \delta + \delta C^{Q+1} \underset{\forall k_{\inmesh{o}} \text{ s.t } k_{\inmesh{o}} \eitherarrow i_{\inmesh{n}}}{\operatorname{mean}} \sum_{j=1}^{L_{Q+1}} |W^d_{k_{\inmesh{o}}j}| \sum_{l_{Q}=1}^{L_{Q}} |W^{(Q)}_{j,l_{Q}}| \cdots \sum_{l_{1}=1}^{L_{1}} |W^{(1)}_{l_{2},l_{1}}| ~ \sum_{j_{\inmesh{o}}=1}^{N_o} |W_{j_{\inmesh{o}}l_1}^e|,\\
&=& \beta + \delta + \delta C^{Q+1} ||\bfmW^{d}||_\infty ||\bfmW^{e}||_\infty \prod_{p=1}^Q ||\bfmW^{(p)}||_\infty.
\neea

\clearpage
\section{Training Details}

\subsection{Network architecture and hyperparameter choices}
\label{app:hyperparams}
$~$

\begin{table}[h]
        \caption{Network architecture and hyperparameter choices for the MOR methods.}
        \label{tab:hyperparams}
    \begin{threeparttable}
        \centering
        \begin{tabular}{c|c|c|c}
             & POD-NN & GCA-ROM & GFN-ROM \\ \hline
             \cellcolor[HTML]{E5E3E3} Bottleneck size ($N$) & \cellcolor[HTML]{E5E3E3} $\lfloor 1.5\times N_\mu \rfloor$ & \cellcolor[HTML]{E5E3E3} $\lfloor 1.5\times N_\mu \rfloor$ & \cellcolor[HTML]{E5E3E3} $\lfloor 1.5\times N_\mu \rfloor$ \\
             Mapper weight ($\omega$) & - & 10 & 10\\
             \cellcolor[HTML]{E5E3E3} Optimiser & \cellcolor[HTML]{E5E3E3} Adam & \cellcolor[HTML]{E5E3E3} Adam & \cellcolor[HTML]{E5E3E3} Adam$^a$\\
             Learning rate & $10^{-3}$ & $10^{-3}$ & $10^{-3}$ \\
             \cellcolor[HTML]{E5E3E3} L2 regularisation & \cellcolor[HTML]{E5E3E3} $10^{-5}$ & \cellcolor[HTML]{E5E3E3} $10^{-5}$ & \cellcolor[HTML]{E5E3E3} $10^{-5}$ \\
             Training epochs & 5000 & 5000 & 5000\\
             \cellcolor[HTML]{E5E3E3} Autoencoder sizes & \cellcolor[HTML]{E5E3E3} - & \cellcolor[HTML]{E5E3E3} [$|\mathcal{M}|$, 200, $N$] & \cellcolor[HTML]{E5E3E3} [$|\mathcal{M}|$, 200, $N$] \\
             Train/test split & 30/70 & 30/70 & 30/70\\
             \cellcolor[HTML]{E5E3E3}Autoencoder activation & \cellcolor[HTML]{E5E3E3} - & \cellcolor[HTML]{E5E3E3} ELU$^b$ & \cellcolor[HTML]{E5E3E3} Tanh$^b$ \\
             Mapper sizes & [$N_\mu$, 50, 50, 50, 50, $N$] & [$N_\mu$, 50, 50, 50, 50, $N$] & [$N_\mu$, 50, 50, 50, 50, $N$] \\
             \cellcolor[HTML]{E5E3E3} Mapper activation & \cellcolor[HTML]{E5E3E3} Tanh & \cellcolor[HTML]{E5E3E3} Tanh$^c$ & \cellcolor[HTML]{E5E3E3} Tanh
        \end{tabular}
        \begin{tablenotes}
          \item[*] POD-G is also undertaken with $N= \lfloor 1.5\times N_\mu \rfloor$ degrees of freedom and a train/test split of 30/70.
          \small
          \item[a] We employ the precomputed adaptive method so that we can use Adam.
          \item[b] Activation not applied to the last decoder layer.
          \item[c] Activation not applied to the last mapper layer.
        \end{tablenotes}
    \end{threeparttable}
\end{table}

\end{document}